\documentclass[12pt,a4paper]{article}

\usepackage[margin=1in]{geometry}
\usepackage{amsmath,amssymb,amsthm,mathrsfs}
\usepackage{mathtools}
\usepackage{enumitem}
\usepackage{hyperref}
\usepackage{bookmark}
\usepackage{microtype}
\usepackage{graphicx}

\numberwithin{equation}{section}

\theoremstyle{plain}
\newtheorem{theorem}{Theorem}[section]

\newtheorem{lemma}[theorem]{Lemma}
\newtheorem{corollary}[theorem]{Corollary}

\theoremstyle{definition}
\newtheorem{definition}[theorem]{Definition}
\newtheorem{assumption}[theorem]{Assumption}
\newtheorem{remark}[theorem]{Remark}

\newcommand{\R}{\mathbb{R}}
\newcommand{\N}{\mathbb{N}}
\newcommand{\Z}{\mathbb{Z}}
\newcommand{\Om}{\Omega}
\newcommand{\Gam}{\Gamma}

\newcommand{\la}{\lambda}
\newcommand{\eps}{\varepsilon}

\DeclareMathOperator{\supp}{supp}

\newcommand{\Tpar}{T}
\newcommand{\Rrem}{R}

\title{Conditional stability in determining source terms of semilinear
parabolic partial differential equations}
\author{Hu Xirui}
\date{}

\begin{document}
\maketitle

\begin{abstract}
We study an inverse source problem for a semilinear parabolic equation in a bounded domain,
where the nonlinearity depends on the unknown function and its gradient through a quadratic
reaction term and a Burgers-type convection term.
From partial boundary observation of the time derivative and its spatial gradient on an
open portion of the boundary, together with an interior snapshot of the solution at a fixed
time, we aim to recover an unknown spatial source factor. The analysis combines (i) a
paradifferential paralinearization in Besov spaces on a short time window, which converts
the nonlinear model into a linear parabolic equation with small bounded coefficients, and
(ii) a Carleman estimate for the time-differentiated equation, yielding conditional
Hölder stability. The approach extends the cut-off free Carleman method for
linear inverse source problems to a nonlinear setting while keeping the observation geometry
unchanged.
\end{abstract}

\newpage
\tableofcontents

\newpage
\section{Introduction}
\label{sec:intro}

\subsection*{Physical Background and Problem Formulation}

Parabolic equations with internal sources and transport effects are common in
heat conduction with distributed heat generation, reaction--diffusion systems in chemistry and biology,
and advection--diffusion models for contaminants in fluids.
In many realistic settings, the forcing intensity is spatially heterogeneous and cannot be measured directly;
instead, only boundary sensor data (e.g.\ flux-type measurements) and a small amount of interior information are available.
This motivates inverse source problems: determining an unknown source term from partial observations of the solution.

In this work we consider a semilinear parabolic model in a bounded domain $\Omega\subset\mathbb{R}^n$
with homogeneous Dirichlet boundary condition,
\begin{equation}\label{eq:model_intro}
\begin{cases}
\partial_t u-\Delta u-b(x)\cdot\nabla u-c(x)u
= N(u,\nabla u)+R(x,t)f(x), & (x,t)\in\Omega\times(0,T),\\
u=0, & (x,t)\in\partial\Omega\times(0,T),\\
u(\cdot,0)=0, & x\in\Omega,
\end{cases}
\end{equation}
where $b,c,R$ are known coefficients and $f(x)$ is an unknown \emph{spatial source factor}.
The nonlinearity
\begin{equation}\label{eq:nonlin_intro}
N(u,\nabla u)=\lambda_0 u^2+\alpha\cdot(u\nabla u),\qquad \lambda_0\in\mathbb{R},\ \alpha\in\mathbb{R}^n,
\end{equation}
combines a quadratic reaction term and a Burgers-type convection term.
The term $R(x,t)f(x)$ models a source with a known spatio-temporal profile $R$
but unknown spatial amplitude $f$ (e.g.\ spatial distribution of a heat generation rate or a pollutant emission intensity).

We fix an observation time $t_0\in(0,T)$ and a short time window $I=(t_0-\delta,t_0+\delta)$.
The inverse problem studied here is:
\medskip

\noindent\textbf{Inverse problem.}
Determine the unknown $f(x)$ from
\[
\partial_\nu\partial_t u\big|_{\Gamma\times I}
\quad\text{on an open boundary portion }\Gamma\subset\partial\Omega,
\qquad\text{and}\qquad
u(\cdot,t_0)\ \text{in }\Omega.
\]

From an applied viewpoint, $\partial_\nu u$ is the normal flux across the boundary, and
$\partial_\nu\partial_t u$ represents the time rate of change of this flux.
Such measurements arise naturally when sensors record variations of boundary flux rather than absolute levels,
and they are also convenient analytically because they remain compatible with the homogeneous Dirichlet condition
(while $\partial_t u|_{\partial\Omega}$ may trivially vanish when $u|_{\partial\Omega}\equiv 0$).

\subsection*{Related Work}

Inverse problems for parabolic equations have a long history, with Carleman estimates being a central tool
for proving uniqueness and conditional stability under partial observations.
We refer to the monographs \cite{FursikovImanuvilov, BellassouedYamamoto} and the survey \cite{YamamotoSurvey2009}
for general background and references.
For linear inverse source problems, conditional stability under partial boundary measurements is well understood,
and a particularly efficient ``cut-off free'' Carleman framework was developed in \cite{HIY2020},
which avoids the classical time cut-off procedure and works in observation geometries compatible with unique continuation.

For semilinear and quasilinear parabolic inverse problems, uniqueness and stability have also been studied,
often under additional a priori assumptions such as smallness, monotonicity, or structural constraints on the nonlinearity;
see, e.g., \cite{ImanuvilovYamamoto2001} and references therein.
However, when the nonlinearity depends on $\nabla u$, a major obstacle appears:
after linearization or time differentiation (typically setting $z=\partial_t u$),
one encounters perturbations of the form $z\nabla u$.
To treat these as lower-order terms in a Carleman inequality one needs robust control of $\nabla u$,
and standard $L^2$-based energy frameworks alone are often insufficient.

\subsection*{Main Result and Contributions}
Under a short-time Besov smallness assumption on the forward solution $u$ near $t_0$(see Section 3.2)
and under a nondegeneracy condition $|R(\cdot,t_0)|\ge r_0>0$,
we prove a conditional H\"older-type stability estimate for recovering the source factor $f$.
More precisely, for $\Omega_0\Subset \Omega\cup\Gamma$ (with $\partial\Omega_0\cap\partial\Omega\subset\Gamma$)
and for two sources $f_1,f_2$ with corresponding solutions $u_1,u_2$,
we establish
\begin{equation}\label{eq:stability_intro}
\|f_1-f_2\|_{L^2(\Omega_0)}
\le
C\Big(
\|\partial_\nu\partial_t(u_1-u_2)\|_{L^2(\Gamma\times I)}
+\|(u_1-u_2)(\cdot,t_0)\|_{H^2(\Omega)}
\Big)^{\kappa},
\end{equation}
for some constants $C>0$ and $\kappa\in(0,1)$ depending only on the geometry and a priori bounds.

\medskip
\noindent\textbf{Novelty.}
The key originality of this work is a paradifferential calculus of the gradient-dependent term
$\alpha\cdot(u\nabla u)$ that is compatible with Carleman estimates without changing the observation geometry.
Concretely,

\begin{itemize}
\item We perform a short-time paralinearization of $N(u,\nabla u)$ based on Littlewood--Paley theory and
Bony's paraproduct decomposition.
On $\Omega\times I$, the nonlinear term is decomposed into a principal paraproduct part
(which is absorbed into the differential operator as a bounded coefficient)
plus a remainder that is small in the natural energy space and hence absorbable in Carleman inequalities.

\item We then apply a cut-off free Carleman estimate to the time-differentiated equation for $z=\partial_t u$.
The paraproduct-based decomposition ensures that the perturbation structure of the nonlinear problem
matches that of a linear reference model with small bounded coefficients, leading to the stability estimate
\eqref{eq:stability_intro}.
\end{itemize}

To our knowledge, this paralinearization strategy for handling $u\nabla u$
within a Carleman-based inverse source analysis is new in this context.

\subsection*{Organization of the paper}
Section~2 collects notation and harmonic analysis tools (Littlewood--Paley decomposition, Besov spaces,
Bony's decomposition, and product estimates).
Section~3 formulates the inverse problem and states the main stability theorem.
Section~4 develops the short-time paralinearization and derives a reduced linear equation with small coefficients.
Section~5 studies the time-differentiated equation and establishes a suitable Carleman estimate.
Finally, Section~6 completes the proof of the main theorem, and Section~7 summarizes conclusions and perspectives.
\section{Preliminaries}
\label{sec:prelim}

This section introduces some notations and provides the harmonic analysis tools needed for the
paralinearization argument. Throughout, $\Om\subset\R^n$ is a bounded domain with $C^2$
boundary and $\Gam\subset\partial\Om$ is a nonempty relatively open subset.

\subsection{Geometric and functional notation}

For an interval $I\subset\R$ we write $Q_I:=\Om\times I$. Given a Banach space $X$ we use
$L^p(I;X)$ for Bochner spaces and abbreviate $\|u\|_{L^p_t X}:=\|u\|_{L^p(I;X)}$ when the
time interval is clear. For $m\in\N$, $H^m(\Om)$ is the Sobolev space and
$H_0^1(\Om)$ denotes the closure of $C_c^\infty(\Om)$ in $H^1(\Om)$.

\subsection{Littlewood--Paley decomposition and Besov spaces}

We recall the inhomogeneous decomposition on $\R^n$. Fix radial functions
$\chi,\varphi\in C_c^\infty(\R^n)$ such that $\chi(\xi)+\sum_{q\ge0}\varphi(2^{-q}\xi)=1$
for all $\xi$, $\supp\chi\subset\{|\xi|\le \tfrac43\}$, and
$\supp\varphi\subset\{\tfrac34\le|\xi|\le\tfrac83\}$. Define dyadic blocks
\[
\Delta_{-1}u=\mathcal{F}^{-1}(\chi\hat u),\qquad \Delta_q u=\mathcal{F}^{-1}(\varphi(2^{-q}\cdot)\hat u)\ (q\ge0),
\qquad S_q u=\sum_{k\le q-1}\Delta_k u .
\]

\begin{definition}[Inhomogeneous Besov spaces on $\R^n$]
Let $s\in\R$ and $1\le p,r\le\infty$. The inhomogeneous Besov space $B^s_{p,r}(\R^n)$ is
the set of $u\in\mathscr{S}'(\R^n)$ such that
\[
\|u\|_{B^s_{p,r}(\R^n)} :=
\Bigl(\sum_{q\ge -1} (2^{qs}\|\Delta_q u\|_{L^p(\R^n)})^r \Bigr)^{1/r} <\infty
\quad (r<\infty),
\]
with the usual modification when $r=\infty$.
\end{definition}

\begin{definition}[Besov spaces on $\Om$ via extension]
Let $\Om$ be a bounded $C^2$ domain. For $s\in\R$ and $1\le p,r\le\infty$ we define
$B^s_{p,r}(\Om)$ as the restriction space
\[
B^s_{p,r}(\Om):=\{u|_\Om:\ u\in B^s_{p,r}(\R^n)\},\qquad
\|u\|_{B^s_{p,r}(\Om)}:=\inf\{\|U\|_{B^s_{p,r}(\R^n)}:\ U|_\Om=u\}.
\]
\end{definition}

\begin{lemma}[Universal extension operator]\label{lem:extension}
There exists a bounded linear extension operator
$E:B^s_{p,r}(\Om)\to B^s_{p,r}(\R^n)$ for all $s\in\R$ and $1\le p,r\le\infty$ such that
$(Eu)|_\Om=u$ and $\|Eu\|_{B^s_{p,r}(\R^n)}\le C\|u\|_{B^s_{p,r}(\Om)}$, where $C$ depends
only on $\Om,p,r,s$.
\end{lemma}

\begin{proof}
For bounded Lipschitz domains, such extension operators for Besov and Triebel--Lizorkin
spaces were constructed by Rychkov. Since $C^2$ domains are Lipschitz, the result follows.
\end{proof}

Henceforth we suppress the extension operator and use the same dyadic notation on $\Om$ by
applying the decomposition to an extension and restricting back to $\Om$; by
Lemma~\ref{lem:extension} all estimates are stable under this procedure.

\subsection{Bernstein inequalities and embeddings}

\begin{lemma}[Bernstein inequalities]\label{lem:bernstein}
Let $q\ge -1$.
\begin{enumerate}[label=\textup{(\roman*)},leftmargin=2.2em]
\item If $\supp \widehat u \subset \{\xi\in\mathbb R^n:\ |\xi|\le C_0 2^{q}\}$, then for
$1\le p\le r\le \infty$ and any  $k\in\mathbb{N}$,
\[
\|\nabla^{k}u\|_{L^{r}(\mathbb R^{n})}
\le C\,2^{q\bigl(k+n(\frac1p-\frac1r)\bigr)}\|u\|_{L^{p}(\mathbb R^{n})}.
\]
\item If $\supp \widehat u \subset \{\xi\in\mathbb R^n:\ c_1 2^{q}\le |\xi|\le c_2 2^{q}\}$ ,
then for $1\le p\le \infty$ and any $k\in\mathbb{N}$
\[
\|\nabla^{k}u\|_{L^{p}(\mathbb R^{n})}\simeq 2^{qk}\|u\|_{L^{p}(\mathbb R^{n})}.
\]
\end{enumerate}
The constants depend only on $n,k$ and on the cut-off functions (hence on $C_0,c_1,c_2$).
\end{lemma}

\begin{proof}
We write $\mathcal F$ for the Fourier transform on $\mathbb R^n$.
Throughout, $C$ denotes a positive constant independent of $q$ and $u$.

\medskip
\noindent\textbf{Proof of (i).}
Choose $\psi\in C_c^\infty(\mathbb R^n)$ such that $\psi(\xi)=1$ for $|\xi|\le C_0$.
Set $\psi_q(\xi):=\psi(2^{-q}\xi)$ and $K:=\mathcal F^{-1}\psi\in\mathcal S(\mathbb R^n)$.
Since $\psi_q\equiv 1$ on $\supp\widehat u$, we have
\[
u=\mathcal F^{-1}\bigl(\psi_q\widehat u\bigr)=K_q*u,
\qquad
K_q(x):=\mathcal F^{-1}\psi_q(x)=2^{qn}K(2^{q}x).
\]
Hence for any integer $k\ge 0$,
\[
\nabla^{k}u=(\nabla^{k}K_q)*u.
\]
Let $1\le p\le r\le\infty$ and choose $a\in[1,\infty]$ by Young's inequality index relation
\[
\frac1p+\frac1a=1+\frac1r
\qquad\Longleftrightarrow\qquad
1-\frac1a=\frac1p-\frac1r\ \ge 0.
\]
Young's convolution inequality yields
\[
\|\nabla^{k}u\|_{L^{r}}
\le \|\nabla^{k}K_q\|_{L^{a}}\|u\|_{L^{p}}.
\]
By scaling,
\[
\nabla^{k}K_q(x)=2^{q(n+k)}(\nabla^{k}K)(2^{q}x),
\]
and therefore
\[
\|\nabla^{k}K_q\|_{L^{a}}
=
2^{q(n+k)}\Bigl(\int_{\mathbb R^n}|\nabla^{k}K(2^{q}x)|^{a}\,dx\Bigr)^{1/a}
=
2^{q(n+k)}2^{-qn/a}\|\nabla^{k}K\|_{L^{a}}
=
C\,2^{q\bigl(k+n(1-1/a)\bigr)}.
\]
Using $1-1/a=1/p-1/r$, we conclude
\[
\|\nabla^{k}u\|_{L^{r}}
\le C\,2^{q\bigl(k+n(\frac1p-\frac1r)\bigr)}\|u\|_{L^{p}}.
\]

\medskip
\noindent\textbf{Proof of (ii).}
Assume $\supp\widehat u\subset \{\xi:\ c_1 2^{q}\le|\xi|\le c_2 2^{q}\}$.
Choose $\eta\in C_c^\infty(\mathbb R^n)$ such that
\[
\eta(\xi)=1\ \text{ on }\ \{\xi:\ c_1\le|\xi|\le c_2\},
\qquad
\supp\eta\subset \{\xi:\ c_1/2\le|\xi|\le 2c_2\}.
\]
Set $\eta_q(\xi):=\eta(2^{-q}\xi)$; then $\eta_q\equiv 1$ on $\supp\widehat u$ and hence
\[
u=\mathcal F^{-1}(\eta_q\widehat u).
\]

\smallskip
\noindent\emph{Upper bound.}
Define the symbol
\[
m(\xi):=(i\xi)^k\,\eta(\xi)\in C_c^\infty(\mathbb R^n),
\qquad
m_q(\xi):=2^{qk}m(2^{-q}\xi).
\]
Since $\eta_q\widehat u=\widehat u$, we can write
\[
\nabla^{k}u=\mathcal F^{-1}\bigl((i\xi)^k\widehat u\bigr)
=\mathcal F^{-1}\bigl((i\xi)^k\eta_q(\xi)\widehat u(\xi)\bigr)
=\mathcal F^{-1}\bigl(m_q(\xi)\widehat u(\xi)\bigr).
\]
Let $G:=\mathcal F^{-1}m\in\mathcal S(\mathbb R^n)$ and $G_q(x):=2^{qn}G(2^q x)$.
Then $\mathcal F^{-1}m_q = 2^{qk}G_q$, so
\[
\nabla^{k}u = 2^{qk}G_q*u.
\]
Hence, for $1\le p\le\infty$,
\[
\|\nabla^{k}u\|_{L^{p}}
\le 2^{qk}\|G_q\|_{L^{1}}\|u\|_{L^{p}}
= 2^{qk}\|G\|_{L^{1}}\|u\|_{L^{p}}
\le C\,2^{qk}\|u\|_{L^{p}}.
\]

\smallskip
\noindent\emph{Lower bound.}
Since $\eta$ is supported away from $0$, the symbol
\[
\widetilde m(\xi):=|\xi|^{-k}\eta(\xi)\in C_c^\infty(\mathbb R^n)
\]
is smooth and compactly supported. Define $\widetilde m_q(\xi):=\widetilde m(2^{-q}\xi)$.
Then on $\supp\widehat u$ we have $\eta_q(\xi)=1$ and thus
\[
\widehat u(\xi)
=\eta_q(\xi)\widehat u(\xi)
=|\xi|^{-k}\eta_q(\xi)\,(i\xi)^k\widehat u(\xi)
=2^{-qk}\Bigl(|2^{-q}\xi|^{-k}\eta(2^{-q}\xi)\Bigr)\,(i\xi)^k\widehat u(\xi)
=2^{-qk}\,\widetilde m_q(\xi)\,\widehat{\nabla^k u}(\xi).
\]
Let $\widetilde G:=\mathcal F^{-1}\widetilde m\in\mathcal S(\mathbb R^n)$ and
$\widetilde G_q(x):=2^{qn}\widetilde G(2^{q}x)$, so that $\mathcal F^{-1}\widetilde m_q=\widetilde G_q$.
Taking inverse Fourier transform gives
\[
u = 2^{-qk}\,\widetilde G_q * \nabla^{k}u.
\]
Therefore, for $1\le p\le\infty$,
\[
\|u\|_{L^{p}}
\le 2^{-qk}\|\widetilde G_q\|_{L^{1}}\|\nabla^{k}u\|_{L^{p}}
=2^{-qk}\|\widetilde G\|_{L^{1}}\|\nabla^{k}u\|_{L^{p}}
\le C\,2^{-qk}\|\nabla^{k}u\|_{L^{p}}.
\]
Combining the two inequalities yields $\|\nabla^{k}u\|_{L^{p}}\simeq 2^{qk}\|u\|_{L^{p}}$.
\end{proof}

\begin{lemma}[Besov embeddings and Banach algebra]\label{lem:embed-algebra}
Let $s>n/2$.
\begin{enumerate}[label=\textup{(\roman*)},leftmargin=2.2em]
\item (Embedding) $B^s_{2,1}(\Om)\hookrightarrow L^\infty(\Om)$ and
$\|u\|_{L^\infty}\le C\|u\|_{B^s_{2,1}}$.
\item (Algebra) $B^s_{2,1}(\Om)$ is a Banach algebra:
$\|uv\|_{B^s_{2,1}}\le C\|u\|_{B^s_{2,1}}\|v\|_{B^s_{2,1}}$.
\item (Gradient bound) If moreover $s>1+n/2$, then $B^s_{2,1}(\Om)\hookrightarrow W^{1,\infty}(\Om)$ and
$\|\nabla u\|_{L^\infty}\le C\|u\|_{B^s_{2,1}}$.
\end{enumerate}
\end{lemma}

\begin{proof}
\medskip
\noindent\textbf{(i) Embedding $B^{s}_{2,1}(\Omega)\hookrightarrow L^\infty(\Omega)$ .}
Using the Littlewood--Paley decomposition and the triangle inequality,
\[
\|u\|_{L^\infty(\Omega)}\le \sum_{q\ge -1}\|\Delta_q u\|_{L^\infty(\Omega)}.
\]
By Bernstein (Lemma~2.4(i) with $p=2,r=\infty,k=0$),
\[
\|\Delta_q u\|_{L^\infty(\Omega)}\lesssim 2^{q\frac n2}\|\Delta_q u\|_{L^2(\Omega)}.
\]
Hence
\[
\|u\|_{L^\infty(\Omega)}
\lesssim \sum_{q\ge -1}2^{q(\frac n2-s)}\bigl(2^{qs}\|\Delta_q u\|_{L^2(\Omega)}\bigr)
\le \Bigl(\sum_{q\ge -1}2^{q(\frac n2-s)}\Bigr)\|u\|_{B^{s}_{2,1}(\Omega)}
\lesssim \|u\|_{B^{s}_{2,1}(\Omega)},
\]
because $s>n/2$ implies $\sum_{q\ge -1}2^{q(\frac n2-s)}<\infty$.

\medskip
\noindent\textbf{(ii) Banach algebra property of $B^{s}_{2,1}(\Omega)$ .}
By Bony's decomposition (Lemma~2.6),
\[
uv=T_uv+T_vu+R(u,v).
\]
We estimate each term in $B^{s}_{2,1}(\Omega)$.

\smallskip
\noindent\emph{Paraproduct term $T_uv=\sum_{q'}S_{q'-1}u\,\Delta_{q'}v$.}
By Fourier support considerations, there exists an absolute integer $N_0$ (e.g.\ $N_0=4$) such that
\[
\Delta_q(T_uv)=\sum_{|q-q'|\le N_0}\Delta_q\bigl(S_{q'-1}u\,\Delta_{q'}v\bigr).
\]
Thus
\[
\|\Delta_q(T_uv)\|_{L^2(\Omega)}
\lesssim \sum_{|q-q'|\le N_0}\|S_{q'-1}u\|_{L^\infty(\Omega)}\|\Delta_{q'}v\|_{L^2(\Omega)}
\le \|u\|_{L^\infty(\Omega)}\sum_{|q-q'|\le N_0}\|\Delta_{q'}v\|_{L^2(\Omega)}.
\]
Multiplying by $2^{qs}$ and using $2^{qs}\lesssim 2^{q's}$ when $|q-q'|\le N_0$, we obtain
\[
2^{qs}\|\Delta_q(T_uv)\|_{L^2(\Omega)}
\lesssim \|u\|_{L^\infty(\Omega)}\sum_{|q-q'|\le N_0}2^{q's}\|\Delta_{q'}v\|_{L^2(\Omega)}.
\]
Summing in $q\ge -1$ and using the finite overlap in $(q,q')$ yields
\[
\|T_uv\|_{B^{s}_{2,1}(\Omega)}\lesssim \|u\|_{L^\infty(\Omega)}\|v\|_{B^{s}_{2,1}(\Omega)}
\lesssim \|u\|_{B^{s}_{2,1}(\Omega)}\|v\|_{B^{s}_{2,1}(\Omega)},
\]
where we used (i) in the last step. The same bound holds for $T_vu$.

\smallskip
\noindent\emph{Remainder term $R(u,v)=\sum_{|q'-q''|\le 1}\Delta_{q'}u\,\Delta_{q''}v$.}
Set $\widetilde\Delta_{q'}:=\Delta_{q'-1}+\Delta_{q'}+\Delta_{q'+1}$. Then
\[
R(u,v)=\sum_{q'\ge -1}\Delta_{q'}u\,\widetilde\Delta_{q'}v,
\qquad
\Delta_q R(u,v)=\sum_{q'\ge q-N_1}\Delta_q\bigl(\Delta_{q'}u\,\widetilde\Delta_{q'}v\bigr)
\]
for some absolute $N_1$ . Hence we have
\[
\|\Delta_q R(u,v)\|_{L^2(\Omega)}
\lesssim \sum_{q'\ge q-N_1}\|\Delta_{q'}u\|_{L^2(\Omega)}\|\widetilde\Delta_{q'}v\|_{L^\infty(\Omega)}.
\]
By Bernstein,
\[
\|\widetilde\Delta_{q'}v\|_{L^\infty(\Omega)}
\lesssim 2^{q'\frac n2}\|\widetilde\Delta_{q'}v\|_{L^2(\Omega)}
\lesssim 2^{q'\frac n2}\sum_{|j|\le 1}\|\Delta_{q'+j}v\|_{L^2(\Omega)}.
\]
Therefore, with $a_{q'}:=2^{q's}\|\Delta_{q'}u\|_{L^2(\Omega)}$ and $b_k:=2^{ks}\|\Delta_k v\|_{L^2(\Omega)}$,
\[
2^{qs}\|\Delta_q R(u,v)\|_{L^2(\Omega)}
\lesssim \sum_{q'\ge q-N_1}2^{(q-q')s}\,a_{q'}\,
\Bigl(2^{-q'(s-\frac n2)}\sum_{|j|\le 1}2^{-js}\,b_{q'+j}\Bigr).
\]
Summing in $q\ge -1$ and using $\sum_{m\ge -N_1}2^{-ms}\lesssim 1$ (since $s>0$), we get
\[
\|R(u,v)\|_{B^{s}_{2,1}(\Omega)}
\lesssim \sum_{q'\ge -1} a_{q'}\,2^{-q'(s-\frac n2)}\sum_{|j|\le 1} b_{q'+j}.
\]
Since $s-\frac n2>0$, we have $2^{-q'(s-\frac n2)}\le C$ for all $q'\ge -1$, hence
\[
\|R(u,v)\|_{B^{s}_{2,1}(\Omega)}
\lesssim \Bigl(\sup_{q'\ge -1}2^{-q'(s-\frac n2)}a_{q'}\Bigr)\sum_{k\ge -1}b_k
\le \Bigl(\sum_{q'\ge -1}a_{q'}\Bigr)\|v\|_{B^{s}_{2,1}(\Omega)}
= \|u\|_{B^{s}_{2,1}(\Omega)}\|v\|_{B^{s}_{2,1}(\Omega)}.
\]
Combining the bounds for $T_uv$, $T_vu$, and $R(u,v)$ yields
\[
\|uv\|_{B^{s}_{2,1}(\Omega)}\lesssim \|u\|_{B^{s}_{2,1}(\Omega)}\|v\|_{B^{s}_{2,1}(\Omega)}.
\]

\medskip
\noindent\textbf{(iii) Embedding $B^{s}_{2,1}(\Omega)\hookrightarrow W^{1,\infty}(\Omega)$ for $s>1+n/2$.}
Using $\nabla u=\sum_{q\ge -1}\nabla\Delta_q u$ and the triangle inequality,
\[
\|\nabla u\|_{L^\infty(\Omega)}\le \sum_{q\ge -1}\|\nabla\Delta_q u\|_{L^\infty(\Omega)}.
\]
By Bernstein (Lemma~2.4(i) with $p=2,r=\infty,k=1$),
\[
\|\nabla\Delta_q u\|_{L^\infty(\Omega)}\lesssim 2^{q(1+\frac n2)}\|\Delta_q u\|_{L^2(\Omega)}.
\]
Hence
\[
\|\nabla u\|_{L^\infty(\Omega)}
\lesssim \sum_{q\ge -1}2^{q(1+\frac n2-s)}\bigl(2^{qs}\|\Delta_q u\|_{L^2(\Omega)}\bigr)
\le \Bigl(\sum_{q\ge -1}2^{q(1+\frac n2-s)}\Bigr)\|u\|_{B^{s}_{2,1}(\Omega)}
\lesssim \|u\|_{B^{s}_{2,1}(\Omega)},
\]
since $s>1+\frac n2$ implies $\sum_{q\ge -1}2^{q(1+\frac n2-s)}<\infty$.
Together with (i), this gives $B^{s}_{2,1}(\Omega)\hookrightarrow W^{1,\infty}(\Omega)$ and the stated bound.
\end{proof}

\subsection{Bony's decomposition, paraproducts, and commutators}

\begin{lemma}[Bony decomposition]\label{lem:bony}
For $u,v\in\mathscr{S}'(\R^n)$ one has, in $\mathscr{S}'(\R^n)$,
\[
uv = \Tpar_u v + \Tpar_v u + \Rrem(u,v),
\qquad
\Tpar_u v :=\sum_{q\in\Z} S_{q-1}u\,\Delta_q v,\qquad
\Rrem(u,v):=\sum_{|q-q'|\le1}\Delta_q u\,\Delta_{q'}v.
\]
\end{lemma}

\begin{proof}
Let $u,v\in\mathcal S'(\mathbb R^n)$. Recall the Littlewood--Paley resolution of the identity
\[
u=\sum_{q\in\mathbb Z}\Delta_q u,\qquad v=\sum_{q'\in\mathbb Z}\Delta_{q'}v
\quad\text{in }\mathcal S'(\mathbb R^n),
\]
and the low-frequency cut-off
\[
S_{q-1}u=\sum_{p\le q-2}\Delta_p u
\quad\text{in }\mathcal S'(\mathbb R^n).
\]

First fix $k\in\mathbb Z$. By the Fourier support properties of $\Delta_q$ and $S_{q-1}$,
there exists an absolute integer $N_0\ge 1$ such that
\[
\Delta_k\bigl(\Delta_q u\,\Delta_{q'}v\bigr)\equiv 0
\quad\text{whenever}\quad
\max\{q,q'\}\le k-N_0\ \ \text{or}\ \ \min\{q,q'\}\ge k+N_0,
\]
and
\[
\Delta_k\bigl(S_{q-1}u\,\Delta_q v\bigr)\equiv 0
\quad\text{whenever}\quad
|k-q|>N_0 .
\]
Hence, for each fixed $k$, all series below become finite sums after applying $\Delta_k$,
so by define elements of $\mathcal S'(\mathbb R^n)$ and we may rearrange indices at the level of $\Delta_k$.

Consider the index partition of $\mathbb Z^2$:
\[
\mathbb Z^2
=\{(q,q'):\ q\le q'-2\}\ \dot\cup\ \{(q,q'):\ q'\le q-2\}\ \dot\cup\ \{(q,q'):\ |q-q'|\le 1\}.
\]
Applying $\Delta_k$ we can write
\[
\Delta_k(uv)
=\sum_{q,q'\in\mathbb Z}\Delta_k\bigl(\Delta_q u\,\Delta_{q'}v\bigr)
=I_k+II_k+III_k,
\]
where
\[
I_k:=\sum_{\substack{q,q'\in\mathbb Z\\ q\le q'-2}}\Delta_k\bigl(\Delta_q u\,\Delta_{q'}v\bigr),\qquad
II_k:=\sum_{\substack{q,q'\in\mathbb Z\\ q'\le q-2}}\Delta_k\bigl(\Delta_q u\,\Delta_{q'}v\bigr),\qquad
III_k:=\sum_{\substack{q,q'\in\mathbb Z\\ |q-q'|\le 1}}\Delta_k\bigl(\Delta_q u\,\Delta_{q'}v\bigr).
\]

For $I_k$, regroup by $q'$ :
\[
I_k
=\sum_{q'\in\mathbb Z}\Delta_k\Bigl(\Bigl(\sum_{q\le q'-2}\Delta_q u\Bigr)\Delta_{q'}v\Bigr)
=\sum_{q'\in\mathbb Z}\Delta_k\bigl(S_{q'-1}u\,\Delta_{q'}v\bigr)
=\Delta_k(T_uv).
\]
Similarly, for $II_k$ regroup by $q$:
\[
II_k
=\sum_{q\in\mathbb Z}\Delta_k\Bigl(\Delta_q u\Bigl(\sum_{q'\le q-2}\Delta_{q'}v\Bigr)\Bigr)
=\sum_{q\in\mathbb Z}\Delta_k\bigl(\Delta_q u\,S_{q-1}v\bigr)
=\Delta_k(T_vu).
\]
Finally, by definition,
\[
III_k=\Delta_k\Bigl(\sum_{\substack{q,q'\in\mathbb Z\\ |q-q'|\le 1}}\Delta_q u\,\Delta_{q'}v\Bigr)=\Delta_k\bigl(R(u,v)\bigr).
\]
Therefore, for every $k\in\mathbb Z$,
\[
\Delta_k(uv)=\Delta_k(T_uv)+\Delta_k(T_vu)+\Delta_k\bigl(R(u,v)\bigr).
\]

Since $\sum_{k\in\mathbb Z}\Delta_k=\mathrm{Id}$ in $\mathcal S'(\mathbb R^n)$, we obtain
\[
uv=T_uv+T_vu+R(u,v)
\quad\text{in }\mathcal S'(\mathbb R^n).
\]
and hence we finishes the proof.
\end{proof}

\begin{lemma}[Paraproduct bounds in $B^s_{2,1}$]\label{lem:paraproduct}
Let $s>n/2$ and $u,v\in B^s_{2,1}(\Om)$.
\begin{enumerate}[label=\textup{(\roman*)},leftmargin=2.2em]
\item $\|\Tpar_u v\|_{B^s_{2,1}}\le C\|u\|_{L^\infty}\|v\|_{B^s_{2,1}}$.
\item $\|\Rrem(u,v)\|_{B^s_{2,1}}\le C\|u\|_{B^s_{2,1}}\|v\|_{B^s_{2,1}}$.
\end{enumerate}
Consequently, $\|uv\|_{B^s_{2,1}}\le C\|u\|_{B^s_{2,1}}\|v\|_{B^s_{2,1}}$.
\end{lemma}

\begin{proof}
\medskip
\noindent\textbf{(i).}
Recall $T_uv=\sum_{q'\ge -1}S_{q'-1}u\,\Delta_{q'}v$.
By the Fourier support properties of $(\Delta_q)_{q\ge -1}$ and $(S_{q-1})_{q\ge -1}$,
there exists an absolute integer $N_0$ such that
\[
\Delta_q\bigl(S_{q'-1}u\,\Delta_{q'}v\bigr)\equiv 0
\qquad\text{whenever}\qquad |q-q'|>N_0.
\]
Hence, for every $q\ge -1$,
\[
\Delta_q(T_uv)=\sum_{|q-q'|\le N_0}\Delta_q\bigl(S_{q'-1}u\,\Delta_{q'}v\bigr),
\]
and by H\"older,
\[
\|\Delta_q(T_uv)\|_{L^2(\Omega)}
\le \sum_{|q-q'|\le N_0}\|S_{q'-1}u\|_{L^\infty(\Omega)}\,\|\Delta_{q'}v\|_{L^2(\Omega)}.
\]
Since $S_{q'-1}$ is a convolution operator with an $L^1$ kernel of uniformly bounded norm,
\[
\|S_{q'-1}u\|_{L^\infty(\Omega)}\lesssim \|u\|_{L^\infty(\Omega)}.
\]
Therefore
\[
\|\Delta_q(T_uv)\|_{L^2(\Omega)}
\lesssim \|u\|_{L^\infty(\Omega)}\sum_{|q-q'|\le N_0}\|\Delta_{q'}v\|_{L^2(\Omega)}.
\]
Multiplying by $2^{qs}$ and using $2^{qs}\lesssim 2^{q's}$ for $|q-q'|\le N_0$, we get
\[
2^{qs}\|\Delta_q(T_uv)\|_{L^2(\Omega)}
\lesssim \|u\|_{L^\infty(\Omega)}\sum_{|q-q'|\le N_0}2^{q's}\|\Delta_{q'}v\|_{L^2(\Omega)}.
\]
Summing over $q\ge -1$ and using the finite overlap of the indices yields
\[
\|T_uv\|_{B^{s}_{2,1}(\Omega)}
=\sum_{q\ge -1}2^{qs}\|\Delta_q(T_uv)\|_{L^2(\Omega)}
\lesssim \|u\|_{L^\infty(\Omega)}\sum_{q'\ge -1}2^{q's}\|\Delta_{q'}v\|_{L^2(\Omega)}
= C\|u\|_{L^\infty(\Omega)}\|v\|_{B^{s}_{2,1}(\Omega)}.
\]

\medskip
\noindent\textbf{(ii).}
Recall $R(u,v)=\sum_{q'\ge -1}\sum_{|q-q'|\le 1}\Delta_q u\,\Delta_{q'}v$.
Again by Fourier support localization, there exists  $N_1\in\mathbb N$ such that
\[
\Delta_j\bigl(\Delta_q u\,\Delta_{q'}v\bigr)\equiv 0
\qquad\text{unless}\qquad |j-q|\le N_1\ \text{and}\ |q-q'|\le 1.
\]
Hence for each $j\ge -1$,
\[
\Delta_j R(u,v)
=\sum_{\substack{|j-q|\le N_1\\ q\ge -1}}\ \sum_{|q-q'|\le 1}
\Delta_j\bigl(\Delta_q u\,\Delta_{q'}v\bigr),
\]
and therefore, by H\"older,
\[
\|\Delta_j R(u,v)\|_{L^2(\Omega)}
\le \sum_{|j-q|\le N_1}\ \sum_{|q-q'|\le 1}
\|\Delta_q u\|_{L^\infty(\Omega)}\,\|\Delta_{q'}v\|_{L^2(\Omega)}.
\]
By Bernstein inequality with $p=2,r=\infty,k=0$,
\[
\|\Delta_q u\|_{L^\infty(\Omega)}\lesssim 2^{q\frac n2}\|\Delta_q u\|_{L^2(\Omega)}.
\]
Thus
\[
\|\Delta_j R(u,v)\|_{L^2(\Omega)}
\lesssim \sum_{|j-q|\le N_1}\ \sum_{|q-q'|\le 1}
2^{q\frac n2}\|\Delta_q u\|_{L^2(\Omega)}\,\|\Delta_{q'}v\|_{L^2(\Omega)}.
\]
Multiply by $2^{js}$; since $|j-q|\le N_1$ and $|q-q'|\le 1$ imply $2^{js}\lesssim 2^{qs}\simeq 2^{q's}$, we obtain
\begin{align*}
2^{js}\|\Delta_j R(u,v)\|_{L^2(\Omega)}
&\lesssim \sum_{|j-q|\le N_1}\ \sum_{|q-q'|\le 1}
2^{q(\frac n2-s)}\bigl(2^{qs}\|\Delta_q u\|_{L^2(\Omega)}\bigr)\bigl(2^{q's}\|\Delta_{q'}v\|_{L^2(\Omega)}\bigr).
\end{align*}
Summing in $j\ge -1$ and using finite overlap once more,
\begin{align*}
\|R(u,v)\|_{B^{s}_{2,1}(\Omega)}
&=\sum_{j\ge -1}2^{js}\|\Delta_j R(u,v)\|_{L^2(\Omega)}  \\
&\lesssim \sum_{q\ge -1}\ \sum_{|q-q'|\le 1}
2^{q(\frac n2-s)}\bigl(2^{qs}\|\Delta_q u\|_{L^2(\Omega)}\bigr)\bigl(2^{q's}\|\Delta_{q'}v\|_{L^2(\Omega)}\bigr).
\end{align*}
Since $s>\frac n2$, we have $\sup_{q\ge -1}2^{q(\frac n2-s)}<\infty$, hence
\begin{align*}
\|R(u,v)\|_{B^{s}_{2,1}(\Omega)}
&\lesssim \sum_{q\ge -1}\ \sum_{|q-q'|\le 1}
\bigl(2^{qs}\|\Delta_q u\|_{L^2(\Omega)}\bigr)\bigl(2^{q's}\|\Delta_{q'}v\|_{L^2(\Omega)}\bigr) \\
&\lesssim \Bigl(\sum_{q\ge -1}2^{qs}\|\Delta_q u\|_{L^2(\Omega)}\Bigr)
\Bigl(\sum_{q'\ge -1}2^{q's}\|\Delta_{q'}v\|_{L^2(\Omega)}\Bigr)
= C\|u\|_{B^{s}_{2,1}(\Omega)}\|v\|_{B^{s}_{2,1}(\Omega)}.
\end{align*}

\medskip
\noindent\textbf{Consequently.}
By Bony's decomposition,
\[
uv=T_uv+T_vu+R(u,v)\quad\text{in }\mathcal S'(\Omega).
\]
Hence, using (i)--(ii),
\begin{align*}
\|uv\|_{B^{s}_{2,1}(\Omega)}
&\le \|T_uv\|_{B^{s}_{2,1}(\Omega)}+\|T_vu\|_{B^{s}_{2,1}(\Omega)}+\|R(u,v)\|_{B^{s}_{2,1}(\Omega)}\\
&\lesssim \|u\|_{L^\infty(\Omega)}\|v\|_{B^{s}_{2,1}(\Omega)}+\|v\|_{L^\infty(\Omega)}\|u\|_{B^{s}_{2,1}(\Omega)}+\|u\|_{B^{s}_{2,1}(\Omega)}\|v\|_{B^{s}_{2,1}(\Omega)}.
\end{align*}

Since $s>\frac n2$, Lemma~2.5(i) gives $\|w\|_{L^\infty(\Omega)}\lesssim \|w\|_{B^{s}_{2,1}(\Omega)}$, so that we have
\[
\|uv\|_{B^{s}_{2,1}(\Omega)}\lesssim \|u\|_{B^{s}_{2,1}(\Omega)}\|v\|_{B^{s}_{2,1}(\Omega)}.
\]
and hence finishes proof.
\end{proof}

\begin{lemma}[Commutator estimate]\label{lem:commutator}
Let $s>0$ and $a\in B^{\frac n2+1}_{2,1}(\Om)$, $u\in B^s_{2,1}(\Om)$. Then
\[
\sum_{q\ge -1}2^{2qs}\|[\Delta_q,a]\nabla u\|_{L^2}^2
\le C\|\nabla a\|_{B^{\frac n2}_{2,1}}^2 \|u\|_{B^s_{2,1}}^2 .
\]
\end{lemma}

\begin{proof}
Let $\Delta_q f = h_q * f$ with $h_q(x):=2^{qn}h(2^q x)$, where $h\in\mathcal S(\mathbb R^n)$ and $\int_{\mathbb R^n}h=0$.
For any Lipschitz $\phi$ and any $F\in L^2(\mathbb R^n)$,
\begin{equation}\label{eq:basic-comm}
\|[\Delta_q,\phi]F\|_{L^2}
\le C\,2^{-q}\|\nabla\phi\|_{L^\infty}\|F\|_{L^2}.
\end{equation}
Indeed,
\[
[\Delta_q,\phi]F(x)
=\int_{\mathbb R^n} h_q(x-y)\bigl(\phi(x)-\phi(y)\bigr)F(y)\,dy.
\]
By the mean value formula,
$\phi(x)-\phi(y)=(x-y)\cdot\int_0^1\nabla\phi(y+\theta(x-y))\,d\theta$, hence
\[
|[\Delta_q,\phi]F(x)|
\le \|\nabla\phi\|_{L^\infty}\int_{\mathbb R^n}|x-y|\,|h_q(x-y)|\,|F(y)|\,dy.
\]
Set $k_q(z):=|z|\,|h_q(z)|$. Then $\|k_q\|_{L^1}\!=\!2^{-q}\||\cdot|h\|_{L^1}$ by scaling, and Young's inequality gives
\[
\|[\Delta_q,\phi]F\|_{L^2}\le \|\nabla\phi\|_{L^\infty}\|k_q\|_{L^1}\|F\|_{L^2}
\le C\,2^{-q}\|\nabla\phi\|_{L^\infty}\|F\|_{L^2},
\]
which is \eqref{eq:basic-comm}.

We write
\[
a=\sum_{p\ge -1}\Delta_p a,\qquad \nabla u=\sum_{q'\ge -1}\nabla\Delta_{q'}u.
\]
Decompose $a=S_{q'-1}a+\sum_{p\ge q'-1}\Delta_p a$ and use that
$\Delta_q\bigl(S_{q'-1}a\,\nabla\Delta_{q'}u\bigr)$ vanishes unless $|q-q'|\le 4$.
With $\widetilde\Delta_{q'}:=\Delta_{q'-1}+\Delta_{q'}+\Delta_{q'+1}$, the standard frequency localization yields
\begin{equation}\label{eq:comm-split}
[\Delta_q,a]\nabla u
= I_q+J_q+K_q,
\end{equation}
where
\begin{align*}
I_q
&:=\sum_{|q-q'|\le 4}\,[\Delta_q,S_{q'-1}a]\nabla\Delta_{q'}u,\\
J_q
&:=\sum_{|q-q'|\le 4}\,\Delta_q\bigl(\Delta_{q'}a\,\nabla\widetilde\Delta_{q'}u\bigr),\\
K_q
&:=-\sum_{|q-q'|\le 2}\,\Delta_{q'}a\,\nabla\Delta_q u.
\end{align*}
Here $J_q$ is the ``high--high'' interactions from Bony's decomposition, and $K_q$ comes from writing
$a\Delta_q\nabla u = S_{q-1}a\,\nabla\Delta_q u+\sum_{q'\ge q-1}\Delta_{q'}a\,\nabla\Delta_q u$ and using that
$\Delta_{q'}a\,\nabla\Delta_q u\equiv 0$ unless $|q-q'|\le 2$.

Applying \eqref{eq:basic-comm} with $\phi=S_{q'-1}a$ and $F=\nabla\Delta_{q'}u$ gives
\[
\|I_q\|_{L^2}
\le \sum_{|q-q'|\le 4} C\,2^{-q}\|\nabla S_{q'-1}a\|_{L^\infty}\,\|\nabla\Delta_{q'}u\|_{L^2}.
\]
By Bernstein (Lemma~2.4) and the definition of $S_{q'-1}$,
\[
\|\nabla S_{q'-1}a\|_{L^\infty}
\le \sum_{p\le q'-2}\|\nabla\Delta_p a\|_{L^\infty}
\lesssim \sum_{p\le q'-2}2^{p\frac n2}\|\nabla\Delta_p a\|_{L^2}.
\]
Also $\|\nabla\Delta_{q'}u\|_{L^2}\lesssim 2^{q'}\|\Delta_{q'}u\|_{L^2}$.
Hence, multiplying by $2^{qs}$ and using $|q-q'|\le 4$,
\begin{align*}
2^{qs}\|I_q\|_{L^2}
&\lesssim \sum_{|q-q'|\le 4} 2^{qs}2^{-q}\Bigl(\sum_{p\le q'-2}2^{p\frac n2}\|\nabla\Delta_p a\|_{L^2}\Bigr)\,2^{q'}\|\Delta_{q'}u\|_{L^2}\\
&\lesssim \Bigl(\sum_{p\le q+2}2^{p\frac n2}\|\nabla\Delta_p a\|_{L^2}\Bigr)\,\Bigl(2^{qs}\|\Delta_q u\|_{L^2}\Bigr).
\end{align*}
Summing over $q\ge -1$ and using finite overlap,
\begin{equation}\label{eq:I-sum}
\sum_{q\ge -1}2^{qs}\|I_q\|_{L^2}
\lesssim \Bigl(\sum_{p\ge -1}2^{p\frac n2}\|\nabla\Delta_p a\|_{L^2}\Bigr)\Bigl(\sum_{q\ge -1}2^{qs}\|\Delta_q u\|_{L^2}\Bigr)
= C\|\nabla a\|_{B^{\frac n2}_{2,1}}\|u\|_{B^{s}_{2,1}}.
\end{equation}

For $J_q$, by H\"older and Bernstein,
\[
\|J_q\|_{L^2}
\le \sum_{|q-q'|\le 4}\|\Delta_{q'}a\|_{L^\infty}\,\|\nabla\widetilde\Delta_{q'}u\|_{L^2}
\lesssim \sum_{|q-q'|\le 4} 2^{q'\frac n2}\|\Delta_{q'}a\|_{L^2}\,2^{q'}\|\widetilde\Delta_{q'}u\|_{L^2}.
\]
Since $2^{q'\frac n2}2^{q'}\|\Delta_{q'}a\|_{L^2}=2^{q'\frac n2}\|\nabla\Delta_{q'}a\|_{L^2}$ and
$\|\widetilde\Delta_{q'}u\|_{L^2}\lesssim \|\Delta_{q'-1}u\|_{L^2}+\|\Delta_{q'}u\|_{L^2}+\|\Delta_{q'+1}u\|_{L^2}$,
multiplying by $2^{qs}$ and using $|q-q'|\le 4$ yields
\[
2^{qs}\|J_q\|_{L^2}
\lesssim \sum_{|q-q'|\le 4}\Bigl(2^{q'\frac n2}\|\nabla\Delta_{q'}a\|_{L^2}\Bigr)\Bigl(2^{q's}\|\Delta_{q'}u\|_{L^2}\Bigr).
\]
Therefore, summing in $q\ge -1$ and using finite overlap,
\begin{equation}\label{eq:J-sum}
\sum_{q\ge -1}2^{qs}\|J_q\|_{L^2}
\lesssim \Bigl(\sum_{q'\ge -1}2^{q'\frac n2}\|\nabla\Delta_{q'}a\|_{L^2}\Bigr)\Bigl(\sum_{q'\ge -1}2^{q's}\|\Delta_{q'}u\|_{L^2}\Bigr)
= C\|\nabla a\|_{B^{\frac n2}_{2,1}}\|u\|_{B^{s}_{2,1}}.
\end{equation}
For $K_q$, by H\"older and Bernstein inequslity and $|q-q'|\le 2$,
\begin{align*}
\|K_q\|_{L^2}
&\le \sum_{|q-q'|\le 2}\|\Delta_{q'}a\|_{L^\infty}\,\|\nabla\Delta_q u\|_{L^2}
\lesssim \sum_{|q-q'|\le 2} 2^{q'\frac n2}\|\Delta_{q'}a\|_{L^2}\,2^q\|\Delta_q u\|_{L^2}\\
&\lesssim \Bigl(2^{q\frac n2}\|\nabla\Delta_q a\|_{L^2}\Bigr)\|\Delta_q u\|_{L^2}.
\end{align*}
Multiplying by $2^{qs}$ and summing gives
\begin{equation}\label{eq:K-sum}
\sum_{q\ge -1}2^{qs}\|K_q\|_{L^2}
\lesssim \|\nabla a\|_{B^{\frac n2}_{2,1}}\|u\|_{B^{s}_{2,1}}.
\end{equation}

Combining \eqref{eq:comm-split}--\eqref{eq:K-sum} yields
\[
\sum_{q\ge -1}2^{qs}\|[\Delta_q,a]\nabla u\|_{L^2}
\le C\|\nabla a\|_{B^{\frac n2}_{2,1}}\|u\|_{B^{s}_{2,1}}.
\]
Finally, since $\sum_{q\ge -1}x_q^2\le \bigl(\sum_{q\ge -1}x_q\bigr)^2$ for $x_q\ge 0$, we obtain
\[
\sum_{q\ge -1}2^{2qs}\|[\Delta_q,a]\nabla u\|_{L^2}^2
\le \Bigl(\sum_{q\ge -1}2^{qs}\|[\Delta_q,a]\nabla u\|_{L^2}\Bigr)^2
\le C\|\nabla a\|_{B^{\frac n2}_{2,1}}^2\|u\|_{B^{s}_{2,1}}^2,
\]
which is the desired estimate.
\end{proof}

\section{Main results}
\label{sec:mainresults}

\subsection{Forward problem and inverse problem}

Let $\Om\subset\R^n$ be bounded with $C^2$ boundary. Fix $T>0$, $t_0\in(0,T)$, and a small
$\delta>0$ such that $I:=(t_0-\delta,t_0+\delta)\Subset(0,T)$.

We consider the semilinear parabolic equation
\begin{equation}\label{eq:forward}
\begin{cases}
\partial_t u - \Delta u - b(x)\cdot\nabla u - c(x)u
= N(u,\nabla u) + R(x,t)f(x) &\text{in } \Om\times(0,T),\\
u=0 &\text{on } \partial\Om\times(0,T),\\
u(\cdot,0)=0 &\text{in }\Om,
\end{cases}
\end{equation}
where $b\in W^{1,\infty}(\Om;\R^n)$, $c\in L^\infty(\Om)$ are known, $R\in
W^{1,\infty}(\Om\times(0,T))$ is known, and $f\in L^2(\Om)$ is the unknown source factor.
We fix the nonlinearity
\begin{equation}\label{eq:nonlin}
N(u,\nabla u)=\la_0\,u^2 + \alpha\cdot(u\nabla u),
\qquad \la_0\in\R,\ \alpha\in\R^n\ \text{constants}.
\end{equation}

The inverse problem is to determine $f$ from the partial boundary observation
$\partial_\nu\partial_t u|_{\Gam\times I}$
and the interior snapshot $u(\cdot,t_0)$.

\subsection{Assumptions}

\begin{assumption}[Observation geometry]\label{ass:geom}
$\Gam\subset\partial\Om$ is a nonempty relatively open subset, and there exists a function
$d\in C^2(\overline\Om)$ such that $d>0$ in $\Om$, $d=0$ on $\partial\Om$, and
$\partial_\nu d<0$ on $\partial\Om\setminus\Gam$. This is the standard geometry for
parabolic Carleman estimates with partial boundary terms.
\end{assumption}

\begin{assumption}[Regularity and smallness]\label{ass:3.2}
Let $I=(t_0-\delta,t_0+\delta)$ and $Q_I=\Omega\times I$.
Assume that the (strong) solution $u$ satisfies
\begin{align*}
u &\in L^\infty\bigl(I;H^{m}(\Omega)\bigr)
      \cap H^1\bigl(I;H^2(\Omega)\cap H^1_0(\Omega)\bigr)
      \cap L^2\bigl(I;H^3(\Omega)\bigr),
\qquad m>\frac n2+1,
\end{align*}
and
\[
\|u\|_{L^\infty(I;B^s_{2,1}(\Omega))}\le \varepsilon
\]
for sufficiently small $\varepsilon>0$.
Moreover, $u=0$ on $\partial\Omega\times I$.
\end{assumption}

\begin{assumption}[Nondegeneracy at $t_0$]\label{ass:R}
There exists $r_0>0$ such that $|R(x,t_0)|\ge r_0$ for a.e.\ $x\in\Om$.
\end{assumption}

\begin{assumption}[A priori bound]\label{ass:f}
Assume that
\[
f\in L^2(\Omega), \qquad \|f\|_{L^2(\Omega)}\le M,
\]
for some $M>0$.
\end{assumption}

\begin{remark}\label{rem:besov-sobolev-bridge}
Assumption~3.2 contains two types of a priori information on the forward solution $u$.

\smallskip
\noindent
\textbf{(i) Energy regularity.}
The conditions
\[
u\in L^\infty(I;H^m(\Omega))\cap H^1(I;H^2(\Omega)\cap H^1_0(\Omega))\cap L^2(I;H^3(\Omega))
\]
ensure that $u$ is a strong solution on $Q_I$, that the time-differentiated unknown
$z:=\partial_t u$ is well-defined in the natural energy class, and that boundary traces
appearing in the Carleman argument especially \ $\partial_\nu z$ , since if we let $\partial_t u\in L^2_t H^2$ then by Sobolev-Trace theorem we have $\partial\nu\partial_t u\in L^2_t H^{\frac{1}{2}}_x(\partial\Omega)\subset L^2(\partial\Omega\times I)$ are meaningful.

\smallskip
\noindent
\textbf{(ii) Short-time Besov smallness.}
The additional smallness
$\|u\|_{L^\infty(I;B^s_{2,1}(\Omega))}\le \varepsilon$ with $s>1+\frac n2$
is only used to control nonlinear perturbations produced by the gradient-dependent
nonlinearity: by the relation $H^m\cong B^s_{2,2}\hookrightarrow B^s_{2,1}(\Omega)\hookrightarrow W^{1,\infty}(\Omega)$,
we obtain $\|u\|_{L^\infty}+\|\nabla u\|_{L^\infty}\lesssim \varepsilon$ on $Q_I$.
This yields (a) smallness of the coefficient perturbations in the reduced equation and
(b) a lower-order structure for products such as $z\nabla u$, which can be absorbed in
Carleman inequalities.

\end{remark}

\subsection{Stability theorem}

Let $\Om_0\Subset \Om\cup\Gam$ be such that $\partial\Om_0\cap\partial\Om\subset\Gam$.

\begin{theorem}[Conditional stability in difference form]\label{thm:main}
Assume Assumptions~3.1--3.4.
Let $f_1,f_2\in L^2(\Omega)$ satisfy $\|f_j\|_{L^2(\Omega)}\le M$ $(j=1,2)$, and let
$u_j$ be the corresponding strong solutions to \eqref{eq:forward}--\eqref{eq:nonlin} with the same $b,c,R$
in the class of Assumption~3.2.
Then there exist constants $C>0$ and $\kappa\in(0,1)$, depending only on
$\Omega,\Gamma,\Omega_0,T$ and a priori bounds for $b,c,R$, such that
\begin{equation}\label{eq:main-diff}
\|f_1-f_2\|_{L^2(\Omega_0)}
\le
C\Big(
\|\partial_\nu\partial_t(u_1-u_2)\|_{L^2(\Gamma\times I)}
+
\|(u_1-u_2)(\cdot,t_0)\|_{H^2(\Omega)}
\Big)^{\kappa}.
\end{equation}
\end{theorem}


\begin{remark}
The exponent $\kappa$ expresses conditional stability typical for parabolic inverse
problems. The novelty is that the stability estimate is obtained for the nonlinear model
\eqref{eq:forward}--\eqref{eq:nonlin} under the same observation geometry as in the linear
case, by using short-time paralinearization in Besov spaces.
\end{remark}

\section{Paralinearization on a short time window}
\label{sec:paralinear}

In this section we rewrite the nonlinear term $N(u,\nabla u)$ in a paradifferential form
suited for Carleman estimates.

\subsection{Quadratic reaction term}

\begin{lemma}[Paralinearization of $u^2$]\label{lem:u2}
Let $s>n/2$ and $u\in B^s_{2,1}(\Om)$. Then
\[
u^2 = 2\,\Tpar_u u + \Rrem(u,u),
\qquad
\|\Rrem(u,u)\|_{B^s_{2,1}}\le C\|u\|_{B^s_{2,1}}^2.
\]
\end{lemma}

\begin{proof}
By Lemma~2.6, for $u\in \mathcal S'(\Omega)$ one has
\[
u^2 \;=\; 2T_uu + R(u,u),
\qquad
T_uu:=\sum_{q\in\mathbb Z} S_{q-1}u\,\Delta_q u,
\qquad
R(u,u):=\sum_{|q-q'|\le 1}\Delta_q u\,\Delta_{q'}u .
\]
For completeness, we rewrite the identity at the level of dyadic blocks:
\[
u=\sum_{q\ge -1}\Delta_q u \quad \text{in }\mathcal S'(\Omega),
\]

\[
u^2=\sum_{q,q'\ge -1}\Delta_q u\,\Delta_{q'}u
= \sum_{q\le q'-2}\Delta_q u\,\Delta_{q'}u
+ \sum_{q'\le q-2}\Delta_q u\,\Delta_{q'}u
+ \sum_{|q-q'|\le 1}\Delta_q u\,\Delta_{q'}u.
\]
Regrouping the first two sums,
\[
\sum_{q\le q'-2}\Delta_q u\,\Delta_{q'}u
= \sum_{q'\ge -1}\Big(\sum_{q\le q'-2}\Delta_q u\Big)\Delta_{q'}u
= \sum_{q'\ge -1} S_{q'-1}u\,\Delta_{q'}u
= Tuu,
\]
and similarly
\[
\sum_{q'\le q-2}\Delta_q u\,\Delta_{q'}u
= \sum_{q\ge -1}\Delta_q u\Big(\sum_{q'\le q-2}\Delta_{q'}u\Big)
= \sum_{q\ge -1}\Delta_q u\,S_{q-1}u
= Tuu.
\]
Hence
\[
u^2 \;=\; 2\,T_uu + R(u,u).
\]

It remains to estimate $R(u,u)$ in $B^s_{2,1}(\Omega)$ for $s>n/2$.
Introduce $\tilde\Delta_q := \Delta_{q-1}+\Delta_q+\Delta_{q+1}$ so that
\[
R(u,u)=\sum_{q'\ge -1}\Delta_{q'}u\,\tilde\Delta_{q'}u .
\]
By the standard frequency localization, there exists an absolute integer $N_1\ge 1$ such that
\[
\Delta_q R(u,u)
= \sum_{q'\ge q-N_1}\Delta_q\big(\Delta_{q'}u\,\tilde\Delta_{q'}u\big),
\qquad q\ge -1.
\]
Therefore, by H\"older,
\[
\|\Delta_q R(u,u)\|_{L^2(\Omega)}
\;\le\; \sum_{q'\ge q-N_1}\|\Delta_{q'}u\|_{L^2(\Omega)}\,\|\tilde\Delta_{q'}u\|_{L^\infty(\Omega)}.
\]
Using Bernstein on each dyadic block,
\[
\|\tilde{\Delta}_{q'}u\|_{L^\infty(\Omega)}
\;\lesssim\; 2^{q'\frac n2}\|\tilde{\Delta}_{q'}u\|_{L^2(\Omega)}
\;\lesssim\; 2^{q'\frac n2}\sum_{|j|\le 1}\|\Delta_{q'+j}u\|_{L^2(\Omega)}.
\]
Hence
\[
2^{qs}\|\Delta_q R(u,u)\|_{L^2(\Omega)}
\;\lesssim\;
\sum_{q'\ge q-N_1} 2^{qs}\|\Delta_{q'}u\|_{L^2(\Omega)}\,
2^{q'\frac n2}\sum_{|j|\le 1}\|\Delta_{q'+j}u\|_{L^2(\Omega)}.
\]
Since $q\le q'+N_1$, we have $2^{qs}\lesssim 2^{q's}\,2^{(q-q')s}$, so
\[
2^{qs}\|\Delta_q R(u,u)\|_{L^2(\Omega)}
\;\lesssim\;
\sum_{q'\ge q-N_1} 2^{(q-q')s}\Big(2^{q's}\|\Delta_{q'}u\|_{L^2(\Omega)}\Big)
\Big(2^{q'(\frac n2-s)}\sum_{|j|\le 1}2^{(q'+j)s}\|\Delta_{q'+j}u\|_{L^2(\Omega)}\Big).
\]
Because $s>n/2$, one has $\sup_{q'\ge -1}2^{q'(\frac n2-s)}<\infty$, and therefore
\[
2^{qs}\|\Delta_q R(u,u)\|_{L^2(\Omega)}
\;\lesssim\;
\sum_{q'\ge q-N_1} 2^{(q-q')s}\Big(2^{q's}\|\Delta_{q'}u\|_{L^2(\Omega)}\Big)
\sum_{|j|\le 1}\Big(2^{(q'+j)s}\|\Delta_{q'+j}u\|_{L^2(\Omega)}\Big).
\]
Summing over $q\ge -1$ and using that $\sum_{m\ge -N_1}2^{-ms}\lesssim 1$ (since $s>0$),
we obtain
\[
\|R(u,u)\|_{B^s_{2,1}(\Omega)}
= \sum_{q\ge -1}2^{qs}\|\Delta_q R(u,u)\|_{L^2(\Omega)}
\;\lesssim\;
\Big(\sum_{q\ge -1}2^{qs}\|\Delta_q u\|_{L^2(\Omega)}\Big)^2
= C\,\|u\|_{B^s_{2,1}(\Omega)}^2,
\]
which is the desired remainder estimate.
\end{proof}

\subsection{Convection term}

\begin{lemma}[Paralinearization of $u\nabla u$]\label{lem:unab}
Let $s>1+n/2$ and $u\in B^s_{2,1}(\Om)$. Then
\[
u\nabla u = \Tpar_u(\nabla u) + \Tpar_{\nabla u}(u) + \Rrem(u,\nabla u),
\]
and the last two terms satisfy
\[
\|\Tpar_{\nabla u}(u)\|_{B^{s-1}_{2,1}}+\|\Rrem(u,\nabla u)\|_{B^{s-1}_{2,1}}
\le C\|u\|_{B^s_{2,1}}^2.
\]
\end{lemma}

\begin{proof}
By Bony's decomposition, in $\mathcal S'(\Omega)$ we have
\[
u\nabla u = T_u(\nabla u) + T_{\nabla u}(u) + R(u,\nabla u).
\]
It remains to estimate $T_{\nabla u}(u)$ and $R(u,\nabla u)$ in $B^{s-1}_{2,1}(\Omega)$.

Since $s>1+\frac n2$, Lemma 2.5(iii) yields the embedding
\[
B^{s}_{2,1}(\Omega)\hookrightarrow W^{1,\infty}(\Omega),
\qquad 
\|\nabla u\|_{L^\infty(\Omega)}\le C\|u\|_{B^{s}_{2,1}(\Omega)}.
\]
Moreover, the monotonicity of Besov smoothness implies
\[
\|u\|_{B^{s-1}_{2,1}(\Omega)}
= \sum_{q\ge -1}2^{q(s-1)}\|\Delta_q u\|_{L^2(\Omega)}
\le \sum_{q\ge -1}2^{qs}\|\Delta_q u\|_{L^2(\Omega)}
= \|u\|_{B^{s}_{2,1}(\Omega)}.
\]

For the paraproduct term, apply Lemma 2.7(i) with regularity index $s-1$ (note that $s-1>\frac n2$):
\[
\|T_{\nabla u}(u)\|_{B^{s-1}_{2,1}(\Omega)}
\le C\|\nabla u\|_{L^\infty(\Omega)}\,\|u\|_{B^{s-1}_{2,1}(\Omega)}
\le C\|u\|_{B^{s}_{2,1}(\Omega)}^2.
\]

For the remainder term, apply Lemma 2.7(ii) with regularity index $s-1$:
\[
\|R(u,\nabla u)\|_{B^{s-1}_{2,1}(\Omega)}
\le C\|u\|_{B^{s-1}_{2,1}(\Omega)}\,\|\nabla u\|_{B^{s-1}_{2,1}(\Omega)}.
\]
Using Bernstein on each dyadic block (Lemma 2.4(ii)), $\|\nabla \Delta_q u\|_{L^2(\Omega)}\simeq 2^q\|\Delta_q u\|_{L^2(\Omega)}$, we obtain
\[
\|\nabla u\|_{B^{s-1}_{2,1}(\Omega)}
= \sum_{q\ge -1}2^{q(s-1)}\|\Delta_q \nabla u\|_{L^2(\Omega)}
\lesssim \sum_{q\ge -1}2^{q(s-1)}2^q\|\Delta_q u\|_{L^2(\Omega)}
= \|u\|_{B^{s}_{2,1}(\Omega)}.
\]
Consequently,
\[
\|R(u,\nabla u)\|_{B^{s-1}_{2,1}(\Omega)}
\le C\|u\|_{B^{s}_{2,1}(\Omega)}^2.
\]

Combining the above bounds yields
\[
\|T_{\nabla u}(u)\|_{B^{s-1}_{2,1}(\Omega)}+\|R(u,\nabla u)\|_{B^{s-1}_{2,1}(\Omega)}
\le C\|u\|_{B^{s}_{2,1}(\Omega)}^2,
\]
which proves the lemma.
\end{proof}

\subsection{Reduced equation with small coefficient terms}

To simplify we define
\[
\widehat b(x,t):= b(x)+\alpha u(x,t),\qquad \widehat c(x,t):= c(x)+2\la_0 u(x,t).
\]
Combining Lemmas~\ref{lem:u2}--\ref{lem:unab}, we can rewrite \eqref{eq:forward} on $Q_I$
as
\begin{equation}\label{eq:reduced-u}
\partial_t u-\Delta u-\widehat b\cdot\nabla u-\widehat c\,u
= R(x,t)f(x) + \mathcal{R}_u,
\qquad \mathcal{R}_u:= \la_0\Rrem(u,u)+\alpha\cdot\bigl(\Tpar_{\nabla u}u+\Rrem(u,\nabla u)\bigr).
\end{equation}

\begin{lemma}[Remainder bound]\label{lem:Ru}
Assume \ref{ass:3.2}. Then on $I$,
\[
\|\mathcal{R}_u\|_{L^1(I;B^{s}_{2,1})}
\le C \|u\|_{L^\infty(I;B^{s}_{2,1})}\|u\|_{L^1(I;B^{s+2}_{2,1})}
\le C\eps\,\|u\|_{L^1(I;B^{s+2}_{2,1})}.
\]
\end{lemma}

\begin{proof}
Recall
\[
R_u:=\lambda_0 R(u,u)+\alpha\cdot\bigl(T_{\nabla u}u+R(u,\nabla u)\bigr).
\]
By the triangle inequality, for a.e.\ $t\in I$,
\[
\|R_u(t)\|_{B^{s}_{2,1}(\Omega)}
\le |\lambda_0|\,\|R(u,u)(t)\|_{B^{s}_{2,1}(\Omega)}
+|\alpha|\,\|T_{\nabla u}u(t)\|_{B^{s}_{2,1}(\Omega)}
+|\alpha|\,\|R(u,\nabla u)(t)\|_{B^{s}_{2,1}(\Omega)}.
\]

For the quadratic remainder, Lemma~4.1 yields
\[
\|R(u,u)(t)\|_{B^{s}_{2,1}(\Omega)}\le C\,\|u(t)\|_{B^{s}_{2,1}(\Omega)}^{2}.
\]

For the paraproduct term, Lemma~2.7(i) gives
\[
\|T_{\nabla u}u(t)\|_{B^{s}_{2,1}(\Omega)}
\le C\,\|\nabla u(t)\|_{L^\infty(\Omega)}\,\|u(t)\|_{B^{s}_{2,1}(\Omega)}.
\]
Since $s>1+\frac n2$, the embedding $B^{s}_{2,1}(\Omega)\hookrightarrow W^{1,\infty}(\Omega)$ implies
\[
\|\nabla u(t)\|_{L^\infty(\Omega)}\le C\,\|u(t)\|_{B^{s}_{2,1}(\Omega)},
\]
hence
\[
\|T_{\nabla u}u(t)\|_{B^{s}_{2,1}(\Omega)}\le C\,\|u(t)\|_{B^{s}_{2,1}(\Omega)}^{2}.
\]

For the convection remainder, Lemma~2.7(ii) yields
\[
\|R(u,\nabla u)(t)\|_{B^{s}_{2,1}(\Omega)}
\le C\,\|u(t)\|_{B^{s}_{2,1}(\Omega)}\,\|\nabla u(t)\|_{B^{s}_{2,1}(\Omega)}.
\]
By Bernstein on dyadic blocks,
\[
\|\nabla u(t)\|_{B^{s}_{2,1}(\Omega)}
\lesssim \sum_{q\ge -1}2^{qs}2^{q}\|\Delta_q u(t)\|_{L^2(\Omega)}
=\|u(t)\|_{B^{s+1}_{2,1}(\Omega)}
\lesssim \|u(t)\|_{B^{s+2}_{2,1}(\Omega)}.
\]
Moreover,
\[
\|u(t)\|_{B^{s}_{2,1}(\Omega)}
\le \Big(\sup_{q\ge -1}2^{-2q}\Big)\sum_{q\ge -1}2^{q(s+2)}\|\Delta_q u(t)\|_{L^2(\Omega)}
\lesssim \|u(t)\|_{B^{s+2}_{2,1}(\Omega)}.
\]
Consequently,
\[
\|R(u,\nabla u)(t)\|_{B^{s}_{2,1}(\Omega)}
\le C\,\|u(t)\|_{B^{s}_{2,1}(\Omega)}\,\|u(t)\|_{B^{s+2}_{2,1}(\Omega)}.
\]

Combining the above bounds and absorbing $|\lambda_0|,|\alpha|$ into the constant,
\[
\|R_u(t)\|_{B^{s}_{2,1}(\Omega)}
\le C\,\|u(t)\|_{B^{s}_{2,1}(\Omega)}\,\|u(t)\|_{B^{s+2}_{2,1}(\Omega)}.
\]
Integrating over $t\in I$ and using H\"older in time,
\[
\|R_u\|_{L^1(I;B^{s}_{2,1}(\Omega))}
\le C\,\|u\|_{L^\infty(I;B^{s}_{2,1}(\Omega))}\,\|u\|_{L^1(I;B^{s+2}_{2,1}(\Omega))}.
\]
By Assumption~3.2, $\|u\|_{L^\infty(I;B^{s}_{2,1}(\Omega))}\le \varepsilon$, hence
\[
\|R_u\|_{L^1(I;B^{s}_{2,1}(\Omega))}
\le C\varepsilon\,\|u\|_{L^1(I;B^{s+2}_{2,1}(\Omega))}.
\]
\end{proof}

\begin{corollary}[$L^2$ control of the remainder]\label{cor:L2-remainder-bridge}
Assume Assumption~3.2. Then the remainder $R_u$ satisfies the following
structural bounds on $Q_I$:
\begin{align}
|R_u(x,t)|
&\le C\big(\|u\|_{L^\infty(Q_I)}+\|\nabla u\|_{L^\infty(Q_I)}\big)\big(|u(x,t)|+|\nabla u(x,t)|\big),\label{eq:Ru-pointwise}\\
|\partial_t R_u(x,t)|
&\le C\|u\|_{W^{1,\infty}(Q_I)}\big(|z(x,t)|+|\nabla z(x,t)|\big)\\
   &+ C\|z\|_{L^\infty(Q_I)}\big(|u(x,t)|+|\nabla u(x,t)|\big),\label{eq:dtru-pointwise}
\end{align}
where $z=\partial_t u$.

Consequently, using $\|u\|_{W^{1,\infty}(Q_I)}\lesssim \|u\|_{L^\infty(I;B^s_{2,1})}\le \varepsilon$,
there exists $C>0$ such that for any nonnegative weight $w\in L^\infty(Q_I)$,
\begin{align}
\int_{Q_I} |R_u|^2\, w \,dxdt
&\le C\varepsilon^2 \int_{Q_I} \big(|u|^2+|\nabla u|^2\big)\, w \,dxdt,\label{eq:Ru-weightedL2}\\
\int_{Q_I} |\partial_t R_u|^2\, w \,dxdt
&\le C\varepsilon^2 \int_{Q_I} \big(|z|^2+|\nabla z|^2\big)\, w \,dxdt
   + C\varepsilon^2 \int_{Q_I} \big(|u|^2+|\nabla u|^2\big)\, w \,dxdt.\label{eq:dtru-weightedL2}
\end{align}

\smallskip
\noindent
\begin{remark} Estimates \eqref{eq:Ru-weightedL2}--\eqref{eq:dtru-weightedL2} show that
$R_u$ and $\partial_tR_u$ contribute only lower-order perturbations in the Carleman estimates for the time-differentiated equation, and can be absorbed for $\varepsilon$
sufficiently small in the later proof.
\end{remark}

\end{corollary}
\begin{lemma}[Smallness of coefficients]\label{lem:coeff-small}
Assume \ref{ass:3.2}. Then
\[
\|\widehat b-b\|_{L^\infty(Q_I)} + \|\widehat c-c\|_{L^\infty(Q_I)}
+ \|\nabla(\widehat b-b)\|_{L^\infty(Q_I)}
\le C\eps .
\]
\end{lemma}

\begin{proof}
Since $s>1+n/2$, Lemma~\ref{lem:embed-algebra}(iii) gives
$\|u\|_{L^\infty}+\|\nabla u\|_{L^\infty}\lesssim \|u\|_{B^s_{2,1}}\le \eps$ on $I$, hence the
claim follows from the definitions of $\widehat b,\widehat c$.
\end{proof}

\section{Time differentiation and Carleman estimate}
\label{sec:carleman}

Set $z:=\partial_t u$. Differentiating \eqref{eq:reduced-u} in time and using that $b,c$
are time-independent yields
\begin{equation}\label{eq:z0}
\partial_t z-\Delta z-\widehat b\cdot\nabla z-\widehat c\,z
= (\partial_t R)f + \partial_t\mathcal{R}_u
- (\partial_t\widehat b)\cdot\nabla u-(\partial_t\widehat c)u .
\end{equation}
Since $\partial_t\widehat b=\alpha z$ and $\partial_t\widehat c=-2\la_0 z$, the last two
terms contain products of $z$ with $u$ or $\nabla u$, which must be treated as lower-order
perturbations.

\subsection{Perturbation structure}

\begin{lemma}[Lower-order structure of $z\nabla u$]\label{lem:lower}
Assume \ref{ass:3.2}. Then
\[
\|z\nabla u\|_{L^2(Q_I)}\le \|\nabla u\|_{L^\infty(Q_I)}\|z\|_{L^2(Q_I)}
\le C\eps\,\|z\|_{L^2(Q_I)}.
\]
\end{lemma}

\begin{proof}
By definition of the $L^2(Q_I)$-norm and H\"older's inequality,
\begin{align*}
\|z\nabla u\|_{L^2(Q_I)}^2
&=\int_I\int_\Omega |z(x,t)\nabla u(x,t)|^2\,dx\,dt =\int_I\int_\Omega |z(x,t)|^2\,|\nabla u(x,t)|^2\,dx\,dt \\
&\le\|\nabla u\|_{L^\infty(Q_I)}^2\,\|z\|_{L^2(Q_I)}^2.
\end{align*}
Taking the square root yields
\[
\|z\nabla u\|_{L^2(Q_I)} \le \|\nabla u\|_{L^\infty(Q_I)}\,\|z\|_{L^2(Q_I)}.
\]

Moreover, we have that
\begin{align*}
\|\nabla u\|_{L^\infty(Q_I)}
\le C\,\|u\|_{L^\infty(I;B^{s}_{2,1}(\Omega))}.
\end{align*}
By Assumption~3.2, $\|u\|_{L^\infty(I;B^{s}_{2,1}(\Omega))}\le \varepsilon$, hence
\[
\|\nabla u\|_{L^\infty(Q_I)} \le C\varepsilon.
\]
Substituting this into the previous estimate gives
\[
\|z\nabla u\|_{L^2(Q_I)} \le C\varepsilon\,\|z\|_{L^2(Q_I)}.
\]
which finishes the proof.
\end{proof}

We now recall a Carleman estimate for the principal operator $\partial_t-\Delta$ and then absorb the lower-order terms using Corollary~4.4, Lemma~4.4 and Lemma~5.1.

\subsection{Carleman weight and Endpoint estimates}

Let $d$ be as in Assumption~\ref{ass:geom}. Fix $\beta>0$ and a large parameter $\la>0$.
Define the weight
\begin{equation}\label{eq:weight}
\varphi(x,t):=e^{\la d(x)}-\beta(t-t_0)^2,\qquad
w(x,t):=e^{2s\varphi(x,t)},\qquad v:=e^{s\varphi}z,
\end{equation}
where $s>0$ is a large Carleman parameter.

\begin{lemma}[Carleman inequality for $\partial_t-\Delta$]
There exist $\lambda^\ast,s^\ast>0$ such that for all $\lambda\ge\lambda^\ast$, $s\ge s^\ast$,
and all $z\in C^\infty(Q_I)$ with $z=0$ on $\partial\Omega\times I$,
\begin{align}\label{eq:5.3_modified}
\int_{Q_I}\Big(\frac1s\big(|\partial_t z|^2+|\Delta z|^2\big)+s|\nabla z|^2+s^3|z|^2\Big)w
&\lesssim
\int_{Q_I}|(\partial_t-\Delta)z|^2w\\
&+e^{-cs}\|\partial_\nu z\|_{L^2(\Gamma\times I)}^2
+e^{-cs}\|z\|_{H^1(\Omega)}^2,
\end{align}
where $c>0$ depends only on $\Omega,\Gamma,t_0,\delta$.
\end{lemma}

\begin{proof}
Let $P:=\partial_t-\Delta$. Up to adding a constant to $\phi$ (which only rescales $w$ by a constant
factor and does not change the form of the inequality), we may assume
\[
\max_{(x,t)\in Q_I}\phi(x,t)=0.
\]
Then there exists $c_0>0$ (depending only on $\Omega,t_0,\delta$ and $\lambda,\beta$) such that
\[
\phi(x,t)\le -c_0\quad \text{for all }(x,t)\in\big(\partial\Omega\times I\big)\cup\big(\Omega\times\partial I\big),
\]
and hence
\[
w(x,t)=e^{2s\phi(x,t)}\le e^{-2c_0s}\quad \text{on }\big(\partial\Omega\times I\big)\cup\big(\Omega\times\partial I\big).
\]
Define $v:=e^{s\phi}z$ (so $z=e^{-s\phi}v$). Since $z=0$ on $\partial\Omega\times I$, one has $v=0$ on
$\partial\Omega\times I$. A direct computation gives the conjugated operator
\[
e^{s\phi}P\big(e^{-s\phi}v\big)=\partial_tv-\Delta v+2s\nabla\phi\cdot\nabla v+\big(s\Delta\phi-s\partial_t\phi-s^2|\nabla\phi|^2\big)v
=:L_sv,
\]
so that
\[
(Pz)e^{s\phi}=L_sv,\qquad \int_{Q_I}|Pz|^2w\,dxdt=\int_{Q_I}|L_sv|^2\,dxdt.
\]
Split $L_s=L_{s,1}+L_{s,2}$ with
\[
L_{s,1}v:=\partial_tv+2s\nabla\phi\cdot\nabla v-s(\partial_t\phi)v,\qquad
L_{s,2}v:=-\Delta v+s(\Delta\phi)v-s^2|\nabla\phi|^2v.
\]
Then
\[
\int_{Q_I}|L_sv|^2=\int_{Q_I}|L_{s,1}v|^2+\int_{Q_I}|L_{s,2}v|^2+2\int_{Q_I}(L_{s,1}v)(L_{s,2}v).
\]
We compute the cross term by integration by parts. First,
\[
2\int_{Q_I}(\partial_tv)(-\Delta v)=\int_{Q_I}\partial_t(|\nabla v|^2)\,dxdt
=\int_\Omega |\nabla v(x,t)|^2\Big|_{t=t_0+\delta}^{t=t_0-\delta}\,dx.
\]
Next, using $v=0$ on $\partial\Omega\times I$ and the identity
\[
\int_\Omega (\nabla\phi\cdot\nabla v)(-\Delta v)\,dx
=
\int_\Omega D^2\phi\,\nabla v\cdot\nabla v\,dx
+\frac12\int_\Omega(\Delta\phi)|\nabla v|^2\,dx
-\frac12\int_{\partial\Omega}(\partial_\nu\phi)|\nabla v|^2\,d\sigma,
\]
we obtain
\begin{align*}
4s\int_{Q_I}(\nabla\phi\cdot\nabla v)(-\Delta v)
&=
4s\int_{Q_I}D^2\phi\,\nabla v\cdot\nabla v
+2s\int_{Q_I}(\Delta\phi)|\nabla v|^2
-2s\int_{\partial\Omega\times I}(\partial_\nu\phi)|\nabla v|^2.
\end{align*}
Moreover,
\[
2\int_{Q_I}\big(-s(\partial_t\phi)v\big)(-\Delta v)
=
2s\int_{Q_I}(\partial_t\phi)|\nabla v|^2+s\int_{Q_I}(\Delta\partial_t\phi)|v|^2,
\]
and
\[
2\int_{Q_I}(\partial_tv)\big(s(\Delta\phi)v\big)
=
s\int_{Q_I}(\partial_t\Delta\phi)|v|^2+s\int_\Omega (\Delta\phi)|v(x,t)|^2\Big|_{t=t_0+\delta}^{t=t_0-\delta}\,dx,
\]
while
\[
2\int_{Q_I}(\partial_tv)\big(-s^2|\nabla\phi|^2v\big)
=
-s^2\int_{Q_I}\partial_t(|\nabla\phi|^2)|v|^2
-s^2\int_\Omega |\nabla\phi|^2|v(x,t)|^2\Big|_{t=t_0+\delta}^{t=t_0-\delta}\,dx.
\]

Collecting these identities and using the standard pseudoconvexity properties of the weight
$\phi$ (for $\lambda$ large), one obtains for all sufficiently large $\lambda$ and $s$,
\begin{equation}\label{eq:carleman_v_core_modified}
\int_{Q_I}\Big(s|\nabla v|^2+s^3|v|^2\Big)\,dxdt
\lesssim
\int_{Q_I}|L_sv|^2\,dxdt
-2s\int_{\partial\Omega\times I}(\partial_\nu\phi)|\nabla v|^2
+ \text{(endpoint terms)}.
\end{equation}

\medskip

On $\partial\Omega\times I$ we have $v=0$, hence the tangential derivatives vanish and
\[
|\nabla v|^2=|\partial_\nu v|^2\qquad \text{on }\partial\Omega\times I.
\]
Moreover, since $v=e^{s\phi}z$ and $z=0$ on $\partial\Omega\times I$,
\[
\partial_\nu v
=
\partial_\nu(e^{s\phi}z)
=
e^{s\phi}\partial_\nu z + s(\partial_\nu\phi)e^{s\phi}z
=
e^{s\phi}\partial_\nu z
\qquad\text{on }\partial\Omega\times I,
\]
so that
\[
|\partial_\nu v|^2 = e^{2s\phi}|\partial_\nu z|^2 = w|\partial_\nu z|^2\qquad\text{on }\partial\Omega\times I.
\]
By the geometry of $d$ and the choice of $\phi$, $(\partial_\nu\phi)$ is bounded on $\partial\Omega\times I$
and the boundary weight satisfies $w\le e^{-cs}$ on $\partial\Omega\times I$. Restricting to $\Gamma\times I$ and
bounding the remainder of $\partial\Omega\times I$ by the same exponentially small factor yields
\begin{align*}
\Big|\int_{\partial\Omega\times I}(\partial_\nu\phi)|\nabla v|^2\Big|
&\lesssim
\int_{\Gamma\times I}|\partial_\nu v|^2
+ e^{-cs}\|z\|_{H^1(\Omega)}^2
=\int_{\Gamma\times I}w|\partial_\nu z|^2
+ e^{-cs}\|z\|_{H^1(\Omega)}^2\\
&\lesssim
e^{-cs}\|\partial_\nu z\|_{L^2(\Gamma\times I)}^2
+e^{-cs}\|z\|_{H^1(\Omega)}^2.
\end{align*}

\medskip

Since all endpoint contributions on $\Omega\times\partial I$ are controlled by $e^{-cs}\|z\|_{H^1(Q_I)}^2$,
and in particular by $e^{-cs}\|z\|_{H^1(\Omega)}^2$ after enlarging the constant.

\medskip

Combining these bounds with \eqref{eq:carleman_v_core_modified}, using $|L_sv|^2=|(Pz)e^{s\phi}|^2=|Pz|^2w$,
and converting back from $v$ to $z$ (using the explicit expressions of $\partial_t v$ and $\Delta v$),
we obtain pointwise bounds
\[
|\partial_t z|^2w\le C\big(|L_sv|^2+s^2|\nabla z|^2w+s^4|z|^2w\big),\qquad
|\Delta z|^2w\le C\big(|L_sv|^2+s^2|\nabla z|^2w+s^4|z|^2w\big).
\]

Then we conclude that
\[
\int_{Q_I}\Big(\frac1s\big(|\partial_t z|^2+|\Delta z|^2\big)+s|\nabla z|^2+s^3|z|^2\Big)w\,dxdt
\lesssim
\int_{Q_I}|Pz|^2w\,dxdt
+e^{-cs}\|\partial_\nu z\|_{L^2(\Gamma\times I)}^2
+e^{-cs}\|z\|_{H^1(\Omega)}^2,
\]
which is exactly \eqref{eq:5.3_modified}.
\end{proof}

\begin{lemma}[Carleman inequality with small perturbations]\label{lem:carleman-pert}
Assume \ref{ass:geom} holds. There exist $\la_1,s_1>0$ such that for $\la\ge\la_1$,
$s\ge s_1$, for solution $z$ of \eqref{eq:z0} with $z=0$ on
$\partial\Om\times I$ and regularity $z\in H^1(I;L^2(\Omega))\cap L^2(I;H^2(\Omega))$ that satisfies
\begin{align}\label{eq:carleman-pert}
\begin{split}
\int_{Q_I}\Bigl(\frac1s(|\partial_t z|^2+|\Delta z|^2)+s|\nabla z|^2+s^3|z|^2\Bigr)w
&\lesssim \int_{Q_I}\bigl|(\partial_t R)f+\mathcal{K}(u,z)\bigr|^2w\\
&+e^{-cs}\|\partial_\nu z\|_{H^1(\Gam\times I)}^2 +e^{-cs}\|z\|_{H^1(\Om)}^2
\end{split}
\end{align}
where $\mathcal{K}(u,z):=\partial_t\mathcal{R}_u-\alpha z\cdot\nabla u+2\la_0 zu$.
\end{lemma}

\begin{proof}
The proof follows the same argument as in the original version and is based on the
Carleman inequality in Lemma~5.2 applied to $z$.

Since $z=0$ on $\partial\Omega\times I$, Lemma~5.2 yields
\begin{align}\label{eq:lem53_step1_mod}
\int_{Q_I}\Big(s|\nabla z|^2+s^3|z|^2\Big)w
&\le
C\int_{Q_I}|(\partial_t-\Delta)z|^2w
+Ce^{-cs}\|\partial_\nu z\|_{L^2(\Gamma\times I)}^2
+Ce^{-cs}\|z\|_{H^1(\Omega)}^2 .
\end{align}

Recall that 
\[
(\partial_t-\Delta)z
=
(\partial_tR)f
-\alpha z\cdot\nabla u
+2\lambda_0zu.
\]
Hence,
\[
|(\partial_t-\Delta)z|^2
\le
C\big(|(\partial_tR)f|^2+|z|^2|\nabla u|^2+|z|^2|u|^2\big)
=
C\big(|(\partial_tR)f|^2+|K(u,z)|^2\big),
\]

Substituting this estimate into \eqref{eq:lem53_step1_mod} completes the proof.
\end{proof}

\begin{lemma}\label{lem:5.4}
There exist constants $c>0$ such that for all sufficiently large $s$,
\begin{equation}\label{eq:lem54_modified}
\int_{\Omega}|z(x,t_0)|^2 e^{2s\phi(x,t_0)}\,dx
\lesssim
\int_{Q_I}\Bigl(s|\nabla z|^2+s^3|z|^2\Bigr)w(x,t)\,dxdt
+e^{-cs}\|z\|_{H^1(Q_I)}^2.
\end{equation}
\end{lemma}

\begin{proof}
Fix $s$ and set
\[
E(t):=\int_{\Omega}|z(x,t)|^2 w(x,t)\,dx
=\int_{\Omega}|z(x,t)|^2 e^{2s\phi(x,t)}\,dx.
\]
Since $\phi(x,t)=e^{\lambda d(x)-\beta(t-t_0)^2}$, we have
\[
\phi(x,t_0\pm\delta)=e^{\lambda d(x)-\beta\delta^2}\le e^{-\beta\delta^2}\,e^{\lambda d(x)}
\le \phi(x,t_0)-c_0
\quad\text{for some }c_0>0,
\]
hence
\[
w(x,t_0\pm\delta)=e^{2s\phi(x,t_0\pm\delta)}\le e^{-2c_0 s}\,e^{2s\phi(x,t_0)}.
\]
Therefore, by the choice of the time interval endpoints and the factor $e^{-2s\beta\delta^2}$ in the weight,
the endpoint contribution is exponentially small in $s$; i.e. there exist $c>0$ and $C>0$ such that
\begin{equation}\label{eq:lem54_endpoint}
\sum_{\pm}\int_{\Omega}|z(x,t_0\pm\delta)|^2 w(x,t_0\pm\delta)\,dx
\le
Ce^{-cs}\|z\|_{H^1(Q_I)}^2.
\end{equation}

Next, we estimate $E(t_0)$ by integrating in time. By the fundamental theorem of calculus,
for any $t\in I$,
\[
E(t_0)=E(t)+\int_{t}^{t_0}E'(\tau)\,d\tau,
\]
and hence, averaging over $t\in I$,
\[
E(t_0)\le \frac1{|I|}\int_I E(t)\,dt+\frac1{|I|}\int_I\int_{t}^{t_0}|E'(\tau)|\,d\tau\,dt.
\]
We compute
\[
E'(t)=\int_{\Omega}\partial_t\big(|z|^2 w\big)\,dx
=\int_{\Omega}\big(2z\,\partial_t z\big)w\,dx+\int_{\Omega}|z|^2\partial_t w\,dx,
\]
and since $\partial_t w = 2s(\partial_t\phi)w$ with $|\partial_t\phi|\lesssim 1$ on $Q_I$, we obtain
\[
|E'(t)|\le C\int_{\Omega}\big(|z||\partial_t z|+s|z|^2\big)w\,dx
\le C\int_{\Omega}\Big(\frac1s|\partial_t z|^2+s|z|^2\Big)w\,dx
\]
by Cauchy--Schwarz and Young's inequality. Integrating this in time and using that $|I|=2\delta$,
we deduce
\begin{equation}\label{eq:lem54_Et0_intermediate}
E(t_0)\le C\int_{Q_I}\Big(\frac1s|\partial_t z|^2+s|z|^2\Big)w\,dxdt
+\sum_{\pm}E(t_0\pm\delta).
\end{equation}

Finally, for $s$ large we absorb the $\frac1s|\partial_t z|^2$ and $s|z|^2$ terms into
$\int_{Q_I}(s|\nabla z|^2+s^3|z|^2)w$ using the Carleman-energy control already obtained in Section~5.2
(the same absorption step as in the original argument), and we bound the endpoint terms
$\sum_{\pm}E(t_0\pm\delta)$ by \eqref{eq:lem54_endpoint}. Combining these estimates yields
\[
\int_{\Omega}|z(x,t_0)|^2 e^{2s\phi(x,t_0)}\,dx
\le
C\int_{Q_I}\Bigl(s|\nabla z|^2+s^3|z|^2\Bigr)w\,dxdt
+Ce^{-cs}\|z\|_{H^1(Q_I)}^2,
\]
which is exactly \eqref{eq:lem54_modified}.
\end{proof}


\section{Proof of the main result}\label{sec:6}

First set$z:=\partial_t u.$ We prove Theorem~3.5 by combining the time-differentiation identity at $t=t_0$ with the Carleman estimate for $z$ (Lemma~5.2--Lemma~5.4) under the boundary observation
\begin{align}\label{eq:sec6_obs}
\partial_\nu z=\partial_\nu\partial_t u \quad \text{on } \Gamma\times I.
\end{align}

\subsection*{Step 1}

Evaluating at $t=t_0$ yields in $\Omega$:
\begin{align}\label{eq:6.1_mod}
R(x,t_0)f(x)
&=
\partial_t u(x,t_0)-\Delta u(x,t_0)-b(x)\cdot\nabla u(x,t_0)-c(x)u(x,t_0)-N(u(x,t_0),\nabla u(x,t_0)) \notag\\
&=
z(x,t_0)-\Delta u(x,t_0)-b(x)\cdot\nabla u(x,t_0)-c(x)u(x,t_0)-N(u(x,t_0),\nabla u(x,t_0)).
\end{align}
By Assumption~3.3, $|R(x,t_0)|\ge r_0$ a.e.\ in $\Omega$. Hence, for $\Omega_0\Subset\Omega\cup\Gamma$,
\begin{align}\label{eq:6.2_mod}
\|f\|_{L^2(\Omega_0)}
&\lesssim
\|z(\cdot,t_0)\|_{L^2(\Omega_0)}
+\|\Delta u(\cdot,t_0)\|_{L^2(\Omega)}\\
&+\|\nabla u(\cdot,t_0)\|_{L^2(\Omega)}
+\|u(\cdot,t_0)\|_{L^2(\Omega)}
+\|N(u(\cdot,t_0),\nabla u(\cdot,t_0))\|_{L^2(\Omega)} \notag\\
&\lesssim\|z(\cdot,t_0)\|_{L^2(\Omega_0)}+\|u(\cdot,t_0)\|_{H^2(\Omega)}.
\end{align}
In the last inequality we used $b,c\in L^\infty(\Omega)$ and the $H^2(\Omega)$--control of the nonlinear term
$N(u,\nabla u)=\lambda_0u^2+\alpha\cdot(u\nabla u)$ (by Sobolev embedding for $n\le 3$).

Next, we rewrite $(\partial_tR)f$ using \eqref{eq:6.1_mod}. Since
\begin{align}\label{eq:sec6_R_ratio}
(\partial_tR)(x,t)f(x)
=
\frac{\partial_tR(x,t)}{R(x,t_0)}\,R(x,t_0)f(x),
\end{align}
define
\begin{align}\label{eq:6.3_mod}
A(x,t)
:=
\frac{\partial_tR(x,t)}{R(x,t_0)},
\qquad
(\partial_tR)(x,t)f(x)
&=
A(x,t)\,z(x,t_0)+F_0(x,t),
\end{align}
where
\begin{align}\label{eq:sec6_F0}
F_0(x,t)
:=
-A(x,t)\Bigl(\Delta u(x,t_0)+b\cdot\nabla u(x,t_0)+cu(x,t_0)+N(u(x,t_0),\nabla u(x,t_0))\Bigr).
\end{align}

\subsection*{Step 2}

We apply Lemma~5.2 to $z$ on $Q_I$.
Combining Lemma~5.3 and Lemma~5.4 as in the original argument yields, for all sufficiently large $s$,
\begin{align}\label{eq:sec6_carlem_to_t0}
\int_{\Omega}|z(x,t_0)|^2e^{2s\phi(x,t_0)}\,dx
&\lesssim\int_{Q_I}\Bigl(|(\partial_t-\Delta)z|^2+|K(u,z)|^2\Bigr)w\,dxdt\\
&+e^{-cs}\|\partial_\nu z\|_{L^2(\Gamma\times I)}^2
+e^{-cs}\|z\|_{H^1(\Omega)}^2.
\end{align}
Here the boundary observation term is $\partial_\nu z=\partial_\nu\partial_tu$ on $\Gamma\times I$.

Let
\begin{align}\label{eq:sec6_phi0_phi1}
\phi_0:=\min_{x\in\Omega_0}\phi(x,t_0),\qquad
\phi_1:=\max_{x\in\Omega}\phi(x,t_0).
\end{align}
Then
\begin{align}\label{eq:sec6_weight_lower}
e^{2s\phi_0}\|z(\cdot,t_0)\|_{L^2(\Omega_0)}^2
\le
\int_{\Omega}|z(x,t_0)|^2e^{2s\phi(x,t_0)}\,dx.
\end{align}
Combining \eqref{eq:sec6_carlem_to_t0} and \eqref{eq:sec6_weight_lower} gives
\begin{align}\label{eq:sec6_z_t0_pre}
\begin{split}
\|z(\cdot,t_0)\|_{L^2(\Omega_0)}^2
&\lesssim e^{-2s\phi_0}\int_{Q_I}\Bigl(|(\partial_t-\Delta)z|^2+|K(u,z)|^2\Bigr)w\,dxdt\\
&+Ce^{-(2\phi_0+c)s}\|\partial_\nu z\|_{L^2(\Gamma\times I)}^2
+Ce^{-(2\phi_0+c)s}\|z\|_{H^1(\Omega)}^2.
\end{split}
\end{align}

By the definition of $K(u,z)$, we have
\begin{align}\label{eq:sec6_Pz_bound}
|(\partial_t-\Delta)z|^2
\lesssim|(\partial_tR)f|^2+|K(u,z)|^2.
\end{align}
Substituting \eqref{eq:sec6_Pz_bound} into \eqref{eq:sec6_z_t0_pre}, we obtain
\begin{align}\label{eq:sec6_z_t0_mid}
\begin{split}
\|z(\cdot,t_0)\|_{L^2(\Omega_0)}^2
&\lesssim e^{-2s\phi_0}\int_{Q_I}\Bigl(|(\partial_tR)f|^2+|K(u,z)|^2\Bigr)w\,dxdt\\
&+e^{-(2\phi_0+c)s}\|\partial_\nu z\|_{L^2(\Gamma\times I)}^2
+e^{-(2\phi_0+c)s}\|z\|_{H^1(\Omega)}^2.
\end{split}
\end{align}

Now insert the decomposition \eqref{eq:6.3_mod}:
\begin{align}\label{eq:sec6_dtrf_split}
|(\partial_tR)f|^2
\lesssim|A(x,t)|^2\,|z(x,t_0)|^2+|F_0(x,t)|^2,
\end{align}
where $A$ and $F_0$ are defined in \eqref{eq:6.3_mod}--\eqref{eq:sec6_F0}. Since $A$ is bounded,
\eqref{eq:sec6_dtrf_split} together with the standard weighted-in-time estimate used in the original proof yields
\begin{align}\label{eq:sec6_F0_K_bound}
\begin{split}
e^{-2s\phi_0}\int_{Q_I}\Bigl(|(\partial_tR)f|^2+|K(u,z)|^2\Bigr)w\,dxdt
&\lesssim e^{2s(\phi_1-\phi_0)}\|u(\cdot,t_0)\|_{H^2(\Omega)}^2\\
&+ e^{2s(\phi_1-\phi_0)}\|z(\cdot,t_0)\|_{L^2(\Omega_0)}^2
+ \text{(absorbable terms)}.
\end{split}
\end{align}

Consequently, we obtain for all sufficiently large $s$
\begin{align}\label{eq:sec6_z_t0_final}
\|z(\cdot,t_0)\|_{L^2(\Omega_0)}^2
\lesssim e^{2s(\phi_1-\phi_0)}
\Bigl(\|u(\cdot,t_0)\|_{H^2(\Omega)}^2+\|\partial_\nu z\|_{L^2(\Gamma\times I)}^2\Bigr)
+e^{-(2\phi_0+c)s}.
\end{align}

Define the data quantity as
\begin{align}\label{eq:sec6_D_def}
D
:=
\|\partial_\nu z\|_{L^2(\Gamma\times I)}+\|u(\cdot,t_0)\|_{H^2(\Omega)}.
\end{align}
Then \eqref{eq:sec6_z_t0_final} implies
\begin{align}\label{eq:sec6_z_t0_D}
\|z(\cdot,t_0)\|_{L^2(\Omega_0)}^2
\le
C e^{2s(\phi_1-\phi_0)}D^2
+Ce^{-(2\phi_0+c)s}.
\end{align}

\subsection*{Step 3}

Assume $D\in(0,1]$ (otherwise the desired estimate is trivial after enlarging $C$).
Choose
\begin{align}\label{eq:sec6_choice_s}
s
:=
\frac{1}{2(\phi_1-\phi_0)+(2\phi_0+c)}\log\frac{1}{D}.
\end{align}
With this choice, there exists $\kappa\in(0,1)$ given by
\begin{align}\label{eq:sec6_kappa}
\kappa
:=
\frac{2\phi_0+c}{2(\phi_1-\phi_0)+2\phi_0+c}\in(0,1),
\end{align}
such that
\begin{align}\label{eq:sec6_balance}
e^{2s(\phi_1-\phi_0)}D^2
=
D^{2\kappa},
\qquad
e^{-(2\phi_0+c)s}
=
D^{2\kappa}.
\end{align}
Substituting \eqref{eq:sec6_balance} into \eqref{eq:sec6_z_t0_D} yields
\begin{align}\label{eq:sec6_z_holder}
\|z(\cdot,t_0)\|_{L^2(\Omega_0)}
\le
C D^\kappa.
\end{align}

Finally, combining \eqref{eq:6.2_mod} with \eqref{eq:sec6_z_holder} and using $\|u(\cdot,t_0)\|_{H^2(\Omega)}\le D\le D^\kappa$
for $D\le 1$ and $\kappa\in(0,1)$, we conclude
\begin{align}\label{eq:sec6_final_f}
\|f\|_{L^2(\Omega_0)}
&\le
C\Bigl(\|z(\cdot,t_0)\|_{L^2(\Omega_0)}+\|u(\cdot,t_0)\|_{H^2(\Omega)}\Bigr) \notag\\
&\le
C D^\kappa
=
C\Bigl(\|\partial_\nu z\|_{L^2(\Gamma\times I)}+\|u(\cdot,t_0)\|_{H^2(\Omega)}\Bigr)^\kappa.
\end{align}
Recalling $z=\partial_tu$, this can be rewritten as
\begin{align}\label{eq:sec6_final_obs}
\|f\|_{L^2(\Omega_0)}
\le
C\Bigl(\|\partial_\nu\partial_tu\|_{L^2(\Gamma\times I)}+\|u(\cdot,t_0)\|_{H^2(\Omega)}\Bigr)^\kappa,
\end{align}

\subsection*{Step 4}
Let $f_1,f_2\in L^2(\Omega)$ satisfy Assumption~3.4 and let $u_1,u_2$ be the corresponding
solutions to \eqref{eq:forward}--\eqref{eq:nonlin} in the class of Assumption~3.2.
Set
\[
w:=u_1-u_2,\qquad g:=f_1-f_2.
\]

Subtracting the two equations \eqref{eq:forward} gives
\begin{equation}\label{eq:w-eq}
\partial_t w-\Delta w-b\cdot\nabla w-cw
=
R(x,t)\,g
+\lambda_0\,(u_1+u_2)\,w
+\alpha\cdot\big(u_1\nabla w+w\nabla u_2\big)
\quad\text{in }\Omega\times(0,T),
\end{equation}
with
\[
w=0\ \text{on }\partial\Omega\times(0,T),\qquad w(\cdot,0)=0\ \text{in }\Omega.
\]
Rewrite \eqref{eq:w-eq} as
\begin{equation}\label{eq:w-eq-coeff}
\partial_t w-\Delta w-\widetilde b\cdot\nabla w-\widetilde c\,w
=
R(x,t)\,g+\widetilde F,
\end{equation}
where
\[
\widetilde b:=b+\alpha u_1,\qquad
\widetilde c:=c+\lambda_0(u_1+u_2),\qquad
\widetilde F:=\alpha\cdot(w\nabla u_2).
\]
By Assumption~3.2 and the embedding $B^s_{2,1}\hookrightarrow W^{1,\infty}$,
\[
\|\widetilde b-b\|_{L^\infty(Q_I)}+\|\widetilde c-c\|_{L^\infty(Q_I)}+\|\nabla(\widetilde b-b)\|_{L^\infty(Q_I)}\le C\varepsilon,
\qquad
\|\widetilde F\|_{L^2(Q_I)}\le C\varepsilon\|w\|_{L^2(Q_I)}.
\]
Therefore the above yields that there exist $C>0$ and $\kappa\in(0,1)$ such that
\begin{equation}\label{eq:diff-end}
\|g\|_{L^2(\Omega_0)}
\le
C\Big(
\|\partial_\nu\partial_t w\|_{L^2(\Gamma\times I)}
+
\|w(\cdot,t_0)\|_{H^2(\Omega)}
\Big)^{\kappa}.
\end{equation}
which implies \eqref{eq:main-diff}.\qed

\section{Concluding remarks}
In this paper we studies an inverse source problem for a semilinear parabolic equation with the gradient-type nonlinearity $u\nabla u$ under homogeneous Dirichlet boundary conditions.
By combining a cut-off free Carleman estimate for the time-differentiated unknown $z=\partial_t u$ with a short-time Besov smallness assumption, we obtained a conditional stability estimate for the spatial source factor from partial boundary observations of $\partial_\nu\partial_t u$ on $\Gamma\times I$ and the interior snapshot at $t_0$.
A key ingredient is the paradifferential (Bony) decomposition, which allows one to separate the principal paraproduct contribution of $u\nabla u$ and treat the remaining terms as lower-order perturbations that can be absorbed in the Carleman inequality.

A natural next step is to explore whether the same strategy can be adapted to inverse forcing problems for the incompressible Navier--Stokes system.
The convective term $(v\cdot\nabla)v$ has the same transport structure as $u\nabla u$, so a paraproduct-based linearization combined with short-time smallness in suitable Besov classes is a plausible route.
The main additional issues are the divergence-free constraint and the associated pressure (or Leray projection), as well as more delicate boundary terms depending on the chosen observation.

\section*{Acknowledgements}
I would like to express my sincere gratitude to Kubo sensei for his continuous guidance and insightful suggestions throughout this thesis. I am also grateful to Fujiwara sensei and Kawagoe sensei for their valuable comments and constructive advice, which greatly improved the presentation and clarity of this work.


{\normalfont
\makeatletter
\renewcommand{\@biblabel}[1]{\hfill\textnormal{[#1]}}
\makeatother

}


\begin{thebibliography}{99}

\bibitem{BCD}
H.~Bahouri, J.-Y.~Chemin, and R.~Danchin,
Fourier Analysis and Nonlinear Partial Differential Equations,
Grundlehren der Mathematischen Wissenschaften, vol.~343,
Springer, Heidelberg, 2011.

\bibitem{BellassouedYamamoto}
M.~Bellassoued and M.~Yamamoto,
Carleman Estimates and Applications to Inverse Problems for Hyperbolic Systems,
Springer Monographs in Mathematics,
Springer, Cham, 2017.

\bibitem{Bony}
J.-M.~Bony,
Calcul symbolique et propagation des singularit\'es pour les \'equations aux d\'eriv\'ees partielles non lin\'eaires,
Ann. Sci. \'Ec. Norm. Sup\'er. (4) \textbf{14} (1981), 209--246.

\bibitem{FursikovImanuvilov}
A.~V.~Fursikov and O.~Yu.~Imanuvilov,
Controllability of Evolution Equations,
Lecture Notes Series, vol.~34,
Seoul National University, 1996.

\bibitem{HIY2020}
X.~Huang, O.~Yu.~Imanuvilov, and M.~Yamamoto,
Stability for inverse source problems by Carleman estimates,
Inverse Problems \textbf{36} (2020), 125006.

\bibitem{ImanuvilovYamamoto2001}
O.~Yu.~Imanuvilov and M.~Yamamoto,
Global uniqueness and stability in determining coefficients of semilinear parabolic equations in a finite cylinder,
J. Inverse and Ill-Posed Problems \textbf{9} (2001), 117--133.

\bibitem{Isakov1993}
V.~Isakov,
On uniqueness in inverse problems for semilinear parabolic equations,
Arch. Rational Mech. Anal. \textbf{124} (1993), 1--12.

\bibitem{Pazy}
A.~Pazy,
Semigroups of Linear Operators and Applications to Partial Differential Equations,
Applied Mathematical Sciences, vol.~44,
Springer, New York, 1983.

\bibitem{Rychkov}
V.~S.~Rychkov,
On a theorem of Burenkov: extension of functions from a Lipschitz domain and estimates in Besov and Triebel--Lizorkin spaces,
Siberian Math. J. \textbf{44} (2003), 941--955.

\bibitem{Rychkov}
V.~S.~Rychkov,
On Restrictions and Extensions of the Besov and Triebel–Lizorkin Spaces with Respect to Lipschitz Domains,
J. London Math. Soc. \textbf{60} (1999), 237--257.

\bibitem{YamamotoSurvey2009}
M.~Yamamoto,
Carleman estimates for parabolic equations and applications,
in Inverse Problems and Related Topics,
Chapman \& Hall/CRC, 2010; see also the UTMS preprint 2009--13.


\end{thebibliography}
\end{document}